\documentclass[12pt, a4paper]{article}
\usepackage[english]{babel}
\usepackage{amsmath, amsthm, amsfonts, amsbsy, amssymb, scrextend, graphicx}
\usepackage{hyperref} 
\usepackage{verbatim}
\usepackage{color}

\usepackage{tikz}
\usepackage[export]{adjustbox}
\usepackage{float}

 \renewcommand{\L}{\mathsf{L}}

\newcommand{\Rr}{\mathsf{R}} 
 \newcommand{\R}{\mathbb{R}}
 \newcommand{\Z}{\mathbb{Z}}
 \newcommand{\Ee}{\mathbb{E}}
\newcommand{\E}{\mathsf{E}}
 
\renewcommand{\P}{\mathsf{P}}
 \newcommand{\Prob}{\mathbb{P}}
 \newcommand{\Probx}{\mathbb{P}_{{\bf x}}}
 \newcommand{\Px}{\mathsf{P}_{{\bf x}}}
 
 \newcommand{\bx}{{\bf x}}

 \newcommand{\ld}{\lambda}
 
\newcommand{\whn}{\widehat n}
\newcommand{\whm}{\widehat m}

\newcommand{\eps}{\varepsilon}
\newcommand{\be}{{\bf e}}

 \newtheorem{thm}{Theorem}
 \newtheorem{lem}{Lemma}
 \newtheorem{cor}{Corollary}
 \newtheorem{rmk}{Remark}
 \newtheorem*{cjt*}{Conjecture}
\newtheorem{prop}{Proposition}
\newtheorem{dfn}{Definition}

\begin{document}

\title{Localisation in a growth model with interaction}

\author{
Marcelo Costa
\footnote{Department of Mathematical Sciences, Durham University, UK. \newline
\indent  Email: m.r.costa@durham.ac.uk
}\\
{\small  Durham University}
\and
Mikhail Menshikov
\footnote{Department of Mathematical Sciences, Durham University, UK. \newline
\indent  Email: mikhail.menshikov@durham.ac.uk
}\\
{\small  Durham University}
\and 
Vadim Shcherbakov
\footnote{Department of Mathematics, Royal Holloway,  University of London, UK. \newline
\indent  Email: vadim.shcherbakov@rhul.ac.uk
}\\
{\small  Royal Holloway,  University of London}
\and
Marina Vachkovskaia
\footnote{Department of Statistics Institute of Mathematics, Statistics and Scientific Computing, University of Campinas, Brazil.
Email: marinav@ime.unicamp.br
}\\
{\small University of Campinas}
}
\maketitle
\begin{abstract}
{\small 
This paper concerns the long term behaviour of  a growth model describing a random sequential allocation 
of  particles  on a finite  cycle graph.
The model can be regarded as a reinforced urn model with graph-based interactions. 
It is motivated by cooperative sequential adsorption, where adsorption rates 
at a site depend on the configuration of existing particles in the neighbourhood of that site.
Our main result is that, with probability one, the growth process will eventually localise  either at a single site, or at a pair of neighbouring sites.
}
\end{abstract}

\section{Introduction}

This paper concerns a  probabilistic model describing a sequential allocation of particles on a finite cycle graph.
The model is motivated by cooperative sequential adsorption (CSA) (see \cite{Ev93}, \cite{Ev97} and references therein).
CSA models are widely applied in physical chemistry for modelling adsorption processes on a material surface onto which particles are  deposited at random.
The main peculiarity of adsorption processes is that deposited particles change adsorption properties of the material.
This motivates  the growth rates defined in equation (\ref{eq:rates}). The growth rates model a particular situation where  
the subsequent particles are more likely to be adsorbed around previously deposited particles.

There is typically a hard-core constraint associated with CSA.
That  is, the adsorption (growth) rate is zero at any location with more than a certain number of  particles.
The asymptotic shape of the spatial configuration of deposited particles  is of primary interest in such models.
Many probabilistic models of spatial growth by monolayer deposition, diffusion and aggregation dynamics present this characteristic. For instance, the Eden model \cite{Ede69}, diffusion-limited aggregation process \cite{WS81}, first-passage percolation models \cite{Ric73} and contact interaction processes \cite{Sch79}.

\begin{figure}[H]
  \centering
  \includegraphics[scale=0.9,trim={5.2cm 20cm 5cm 5cm},clip]{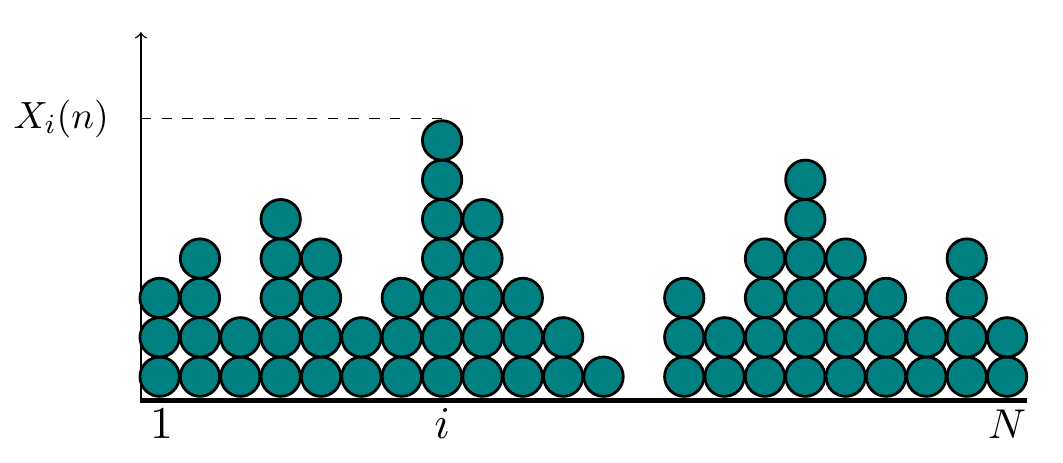}
  \vspace*{-1mm}
  \caption{\small{Multilayer adsorption/random deposition model}}
  
\label{fig1}
\end{figure}

In contrast, in our model (defined in Section \ref{main}) we allow any number of particles to be deposited at each site.
This is motivated by growing interfaces (Figure \ref{fig1}) associated with multilayer adsorption processes (see \cite{BP91}, \cite{JM87} and \cite{Pen08-2}).
Even though the random nature of these processes is usually emphasized in the physical literature, there is a limited number of rigorous formulations and published results in this field (most of them in \cite{Pen08} and \cite{PY02}).
Our model is closely related to a variant of random deposition models, but as we do not apply any of the techniques from this field, we refer the reader to the survey on surface growth \cite{BS95}.

Our  model can be naturally  interpreted in terms of  interacting urn models.
In the case of no interaction, in which the growth rate at site $i$ is given by $\Gamma(x_i)$, where $x_i$ is the number of existing particles at site $i$ and $\Gamma : \Z_+ \rightarrow (0, \infty)$ is a given function (called the reinforcement rule \cite{CCL13} or feedback function \cite{RO07}), our model coincides with a generalised P\'olya urn (GPU)  model with a particular reinforcement rule $\Gamma$.
 Each site (with no underlying graph structure) corresponds to a different colour of ball. The growth rule corresponds to choosing an existing ball of colour $i$, with probability proportional to $\Gamma(x_i)$, and adding a new ball of that colour. The case $\Gamma(x) = x$ is the classical P\'olya urn.

Many non-interacting urn models' techniques rely on the so called Rubin's exponential embedding (first appearing in \cite{Dav90}), which classifies the two possible limiting behaviours.
 First, there almost surely exists a site $i$ that gets all but finitely many particles. Second, the number of particles at every site grows almost surely to infinity.
For a comprehensive survey on GPU models, variants and applications, see \cite{Pem07} and references therein. 

\begin{figure}[H]
  \centering
  \includegraphics[scale=0.8, trim={4cm 19cm 4cm 4cm},clip]{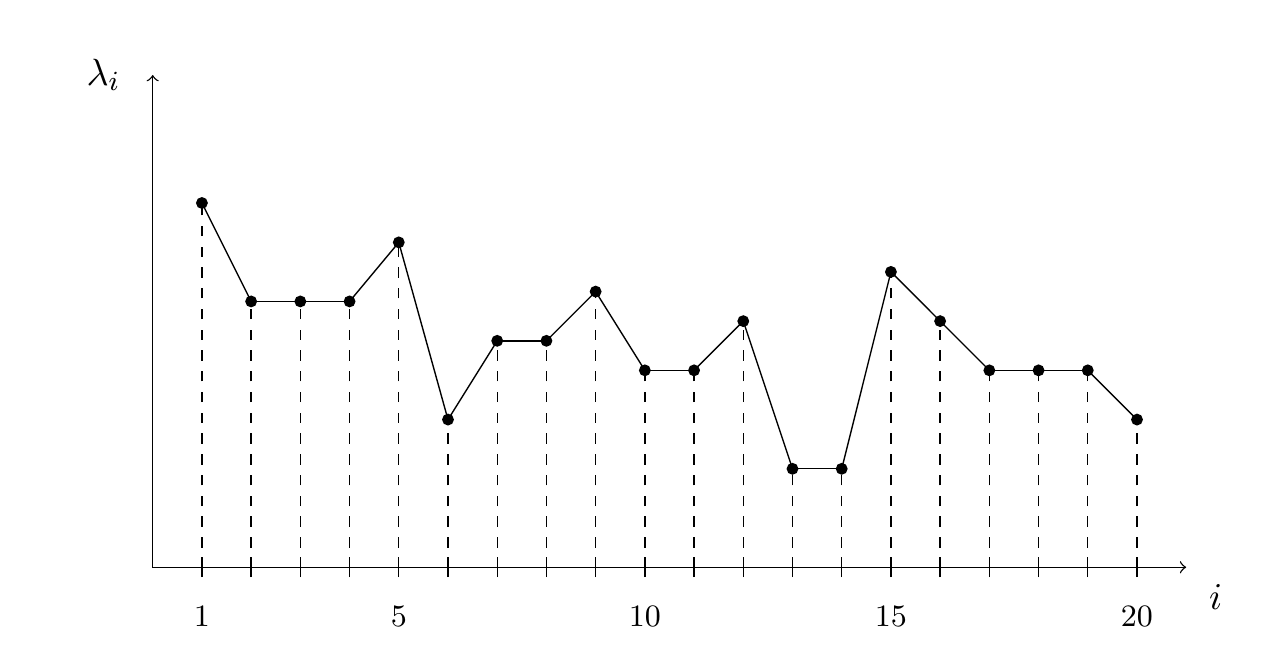}
  \vspace*{-4mm}
  \caption{\small{An interpolated graph of a particular parameter set $(\lambda_i)_{i=1}^{20}$.}}
  \label{fig2}
\end{figure}

In contrast, we consider growth rules with graph-based interactions (as in \cite{SV10}) where the underlying graph is a cycle with $N$ sites.
In our growth model the rate of growth at site $i$ is given by a site-dependent reinforcement rule $\Gamma_i = \exp(\lambda_i u_i)$, where $\lambda_i >0$  and $u_i$ is the number of existing particles in a neighbourhood of site $i$.
This allows one to take into account the case where different sites might possibly have different reinforcement schemes (Figure \ref{fig2}). In other words, the case where each site has its own intrinsic `capacity' parameter, which is what would be expected in many real-life situations.
Although the model can easily be defined for a general graph, the results will heavily depend on its topological properties. In this paper we only address the case of a cycle graph.
See \cite{BBCL15} and \cite{FFK11} for results on general graphs but different growth rules.
The main result of the present paper classifies, in terms of the set of parameters $\Lambda = (\lambda_i)_{i=1}^N$, the two possible behaviours of the model. The first behaviour is localization of growth at a single site.
This means that from a random moment of time onwards, all subsequent particles are allocated at a particular site.
The second is localization of growth at a pair of neighbouring sites with equal $\lambda$ parameter.
Similarly as in the first case, this means that from a random moment of time onwards, all subsequent particles are allocated at a particular pair of neighbouring sites.
In particular, if $\lambda_i\neq \lambda_{i+1}$ for all $i$, then, with probability one, the growth will eventually localise at a single site.
On the other hand, if $\lambda_i\equiv \lambda$, then, with probability one, the growth will eventually localise at a pair of neighbouring sites. 
In the general case of a fixed and arbitrary parameter set $\Lambda$, only the above two types of limiting behaviour  are possible.
Theorem~\ref{T1} below provides a complete characterization of the parameter set $\Lambda$ and associated subsets where only one of the regimes, or both, may happen.  

The model with $\Gamma_i = \exp(\lambda u_i)$, i.e. $\lambda_i\equiv \lambda  \in \mathbb{R}$, was first considered in \cite{SV10}, and an analogue  of Theorem \ref{T1} (Theorem 3 in~\cite{SV10}) was proved for this particular case of site-independent  parameter $\lambda$.

The paper is organised as follows. In Section \ref{main}, we formally define the model, fix some terminology and state Theorem~\ref{T1} which is our main result. 
The proof of the theorem appears in Section \ref{proofT1} and relies essentially on Lemmas~\ref{L1}-\ref{L8} stated in Section \ref{lemmas} and proved in Section~\ref{ProofsLemmas}.  
Section \ref{App} contains results concerning sums of random geometric progressions, which are  of interest in their own right.
These results combined with stochastic domination techniques are constantly used in the proofs of Lemmas~\ref{L5}-\ref{L8}.

\section{The model and main result}
\label{main}

Consider a cycle graph with $N \geq 4$ vertices (sites)
 enumerated by the first $N$ natural numbers such that  $1\sim 2 \sim \ldots \sim N-1 \sim N\sim1$, where $i\sim j$  indicates that sites $i$ and $j$ are incident. 
Let $\Z_{+}$ be the set of non-negative integers and $\Lambda=\{\lambda_1,...,\lambda_N\}$ be an arbitrary set of positive real numbers. 
Given $\bx=(x_1, \ldots, x_N)\in \Z_{+}^N$, define the growth rates as
\begin{equation}
\label{eq:rates}
    \Gamma_i({\bf x})= e^{\ld_i(x_i+\sum_{j\sim i}x_j)},\quad i=1,\ldots, N.
\end{equation}
Consider a discrete-time Markov chain $X(n)=(X_1(n),\ldots, X_N(n))\in\Z_{+}^N$ with the following 
transition probabilities
\begin{equation*}
\P(X_i(n+1)=X_i(n)+1|X(n)=\bx)=\frac{\Gamma_{i}(\bx)}{\sum_{k=1}^N\Gamma_k(\bx)},\, i=1,\ldots, N,\, \bx\in \Z_{+}^N.
\end{equation*}
The Markov chain describes the evolution of the number of particles sequentially allocated at each site of the graph.
Given the configuration of particles  $X(n)=\bx \in \Z^N_+$  at time $n$,
the next incoming particle is placed at site $i$ with probability proportional to   $\Gamma_i({\bf x})$.

\begin{dfn}
\label{D1}
For $i\in\{1,\ldots, N\}$ and auxiliary indices $modulo \>  N$
\begin{enumerate}
\item a site $\{i\}$ is a local minimum, if $\lambda_{i}<\min(\lambda_{i-1}, \lambda_{i+1})$;
\item a pair of sites $\{i, i+1\}$ is  a local minimum of size $2$, if
$\lambda_i =\lambda_{i+1}<\min(\lambda_{i-1}, \lambda_{i+2})$;
\item  a site $\{i\}$ is a local maximum,  if $\lambda_{i}>\max(\lambda_{i-1}, \lambda_{i+1})$;
\item a pair  of sites $\{i, i+1\}$ is   a saddle point, if 
$$\min(\lambda_{i-1}, \lambda_{i+2})<\lambda_i=\lambda_{i+1}<\max(\lambda_{i-1}, \lambda_{i+2});$$
\item a site $\{i\}$ is  a growth  point, if either $\lambda_{i-1}<\lambda_{i}<\lambda_{i+1}$, 
or $\lambda_{i-1}>\lambda_{i}>\lambda_{i+1}$.
\end{enumerate}
\end{dfn}

\begin{dfn}
\label{D2}
Let  $\{i, i+1\}$ be a local minimum of size two. We say that it is a local minimum of size 2 and 
\begin{enumerate}
\item[ 1)]    type $1$, if 
$\lambda_i=\lambda_{i+1}> \frac{\lambda_{i-1}\lambda_{i+2}}{\lambda_{i-1}+\lambda_{i+2}},$
\item[ 2)] 
type $2$, if 
$\lambda_i=\lambda_{i+1}\leq \frac{\lambda_{i-1}\lambda_{i+2}}{\lambda_{i-1}+\lambda_{i+2}}$.
\end{enumerate}
\end{dfn}

The following theorem is the main result of the paper.

\begin{thm}
\label{T1}
 For every $\bx \in \Z_+^N$ and
\begin{itemize}
\item[ i)] for every local maximum $\{k\}$, 
with positive probability, 
$$\lim_{n\to \infty} X_i(n)=\infty \text{ if and only if } i=k;$$
\item[ ii)] for every pair $\{k, k+1\}$ where $\lambda_k = \lambda_{k+1}=:\lambda$,
 but not a local minimum of size 2 and type 2, with positive probability, 
\begin{align*}
&\lim_{n\to \infty} X_i(n)=\infty,\text{ if and only if }i\in \{k, k+1\},\text{ and} \\
&\lim_{n\to\infty} \frac{X_{k+1}(n)}{X_{k}(n)}=e^{\lambda R},
\end{align*}
 where $R=\lim_{n\to\infty} [X_{k+2}(n)-X_{k-1}(n)]\in \Z$.
\end{itemize}
No other limiting behaviour is possible. That is, 
with probability 1, exactly one of the above events occurs in a random location $\{k\}$ or $\{k,k+1\}$ as described in i)  and ii), respectively.
\end{thm}

\section{Lemmas}
\label{lemmas}

We start  with notations that will be used throughout the proofs.
Given  $i = 1,..., N$, define the following events
\begin{align*}
 A_n^i& := \{\text{at time } n  \text{ a particle is placed  at site } i\}, \,  n \in \Z_+, \\
 A_n^{i, i+1}&:=  A_n^i \cup  A_n^{i+1}, \,  n \in \Z_+.
\end{align*}
Define also the following events
\begin{align*}
A_{[n_1,n_2]}^i &:= \bigcap_{n= n_1}^{n_2} A_n^i,\\
A^{i, i+1}_{[n_1, n_2]}&:=\bigcap_{n=n_1}^{n_2}A^{i, i+1}_n,
\end{align*}
indicating that from time $n_1$ to $n_2$ all particles are placed at site $i$, and at sites $i$ or $i+1$, respectively. 
Further, events $A_{[n, \infty)}^i$ and $A_{[n, \infty)}^{i, i+1}$ denote the corresponding limiting cases as $n_2$ goes to infinity.

Let ${\bf e}_i\in \Z_{+}^N$ be a vector, whose $i$-th coordinate is $1$, and all other coordinates are zero.
Given $\bx \in \Z_+^N$, define the following  probability measure 
$\Px(\cdot) = \P(\, \cdot \, | \, X(0)={\bf x})$.

In lemmas and proofs below we denote by $\epsilon$ and $\varepsilon$, possibly with subscripts, various positive constants whose values depend only on $N$ and $(\lambda_i)_{i=1}^ N$ and may vary from line to line. Also, the results are stated only for the essentially different cases, and whenever there are trivially symmetric  situations (e.g. $\lambda_{k-1} < \lambda_k < \lambda_{k+1}$ and $\lambda_{k-1} > \lambda_k > \lambda_{k+1}$), we state and prove only one of them in order to avoid unnecessary repetition.

\begin{lem}
\label{L1} 
Suppose that  $\{k\}$ is a local maximum, and $\bx\in \Z_{+}^N$ is such that 
$\Gamma_k(\bx)=\max_{i}\Gamma_i(\bx)$. Then, with positive probability, all subsequent particles are allocated at $k$, 
i.e. $\Px\left(A_{[1,\infty)}^k\right)\geq \epsilon$ for some  $\epsilon>0$.
\end{lem}
Lemma \ref{L1} describes the only case where localisation of growth at a single site can occur, namely, at a local maximum. 
\begin{lem}
\label{L2} 
Suppose that  $\{k\}$ is a growth  point,  and $\bx\in \Z_{+}^N$ is such that 
$\Gamma_k(\bx)=\max_{i}\Gamma_i(\bx)$.
 If $\lambda_{k-1}<\lambda_k<\lambda_{k+1}$, then there exist $n = n(\bx, \Lambda)\in\Z_{+}$ and $\epsilon> 0$,
 such that 
$\Px\left(A_{[1,n]}^k\right) \geq \epsilon$
and $\Gamma_{k+1}(\bx+n\be_k)=\max_{i}\Gamma_i(\bx+n\be_k)$.
\end{lem}

\begin{lem}
\label{L3}
Suppose that  $\{k\}$ is  a local minimum, and  $\bx\in \Z_{+}^N$ is such that  $\Gamma_k(\bx)= \max_{i}\Gamma_i(\bx)$.
 Then there exist   $n=n(\bx, \Lambda)\in \Z_{+}$ and  $\epsilon>0$, 
 such that  $\Px\left(A_{[1,n]}^k\right) \geq \epsilon$ and 
$\max(\Gamma_{k-1}(\bx+n\be_k), \Gamma_{k+1}(\bx+n\be_k))=\max_{i}\Gamma_i(\bx+n\be_k).$
\end{lem}
Lemmas~\ref{L2}-\ref{L3} describe the following effect.
If the maximal rate is attained at a site which is either a growth point or a local minimum, then, with positive probability, allocating $n=n(\bx, \Lambda)$ particles at that site results in relocation  of the maximal rate to a nearest neighbour with larger parameter $\lambda$.
It should be noted that the number of particles required for relocation (the relocation time) is deterministic and depends only on the starting configuration $\bx$ and parameter set $\Lambda$. 

\begin{lem}
\label{L4}
Suppose that   
$\Gamma_{k}(\bx)=\max_{i}\Gamma_i(\bx)$.
\begin{itemize}
\item[ 1)] $\lambda_{k-1}<\lambda_k=\lambda_{k+1} \geq \lambda_{k+2}$; or
\item[ 2)]  $\lambda_{k-1}=\lambda_k=\lambda_{k+1}\geq \lambda_{k+2},$ and 
 $\Gamma_{k+1}(\bx)\geq \Gamma_{k-1}(\bx)$, 
\end{itemize}
then, with positive probability, all subsequent particles are allocated at sites $\{k, k+1\}$, i.e. 
$\Px\left(A_{[1,\infty)}^{k, k+1}\right) \geq \epsilon$ for some  $\epsilon>0$.
\end{lem}

 Lemma~\ref{L4} describes a case of the second possible limiting behaviour of the model, i.e. localisation of growth
at a pair of neighbouring sites. 

\begin{dfn}
\label{D3}
Define the following stopping times
\begin{align*}
\tau_k&=\inf(n: X_k(n)=X_k(0)+1),\\  
w_{k}^+&=\min(\tau_i :  \,\, i\neq k, k+1, k+2),\\
w_{k}^-&=\min(\tau_i :  \,\, i\neq k-1, k, k+1),  \text{ for } k=1, \ldots, N,
\end{align*}
where the usual convention is that $$\inf(\emptyset)=\infty\text{ and } \min(a, \infty)=a,
\text{ for } a\in\R_{+}\cup \{\infty\}.$$
\end{dfn}

The above stopping times and the quantities $r, z_1$ and $z_2$ below will appear throughout Lemmas \ref{L4}-\ref{L8} and their proofs. 

\begin{dfn}
\label{D4}
Given $\bx \in \Z_+^N$ define 
\begin{equation}
\label{r}
r:=r(\bx)=x_{k+2}-x_{k-1}.
\end{equation}
In addition, if a pair of sites $\{k, k+1\}$ is such that $\lambda_k=\lambda_{k+1}=:\lambda$
and
\begin{align}
\label{z}
  &\lambda_{k-1}>\lambda, \text{ define  } z_1=\frac{1}{\lambda}\log\left(\frac{\lambda_{k-1}-\lambda}{\lambda}\right),\\
  &\lambda_{k+2}>\lambda, \text{ define  } z_2=\frac{1}{\lambda}\log\left(\frac{\lambda}{\lambda_{k+2}-\lambda}\right).
\end{align}
\end{dfn}

Before stating Lemma \ref{L5}, let us denote by $B_k$ the event in which a particle arrives in finite time at $k+2$ before anywhere 
outside $\{k,k+1,k+2\}$. That is to say, 
\begin{equation}
\label{Bk}
B_k:=\{\tau_{k+2}<w_{k}^+\}.
\end{equation}

\begin{lem}
\label{L5}
Suppose that  a pair of sites  $\{k, k+1\}$ is a saddle point with
$\lambda_{k-1}< \lambda_k=\lambda_{k+1}=:\lambda<\lambda_{k+2}$, and 
$\bx\in  \Z_{+}^N$ is such that 
\begin{equation}
\label{A}
\max(\Gamma_{k}(\bx), \Gamma_{k+1}(\bx))=\max_{i}\Gamma_i(\bx).
\end{equation}
\begin{enumerate}
\item[1)] Then there exists  $\epsilon>0$ such that  
$$
\Px\left(A_{[1,\infty)}^{k, k+1}\bigcup B_k\right)
=\Px\left(A_{[1,\infty)}^{k, k+1}\right)
+\Px\left( B_k\right)\geq \epsilon.
$$
\item[2)] If $r<z_2$, 
then,  with positive probability, all subsequent particles are allocated at sites $\{k, k+1\}$, 
i.e. $\Px\left(A_{[1,\infty)}^{k, k+1}\right)\geq \varepsilon$ for some   $\varepsilon>0$.

\item[ 3)] If  
$r\geq z_2$,
then  $\Px\left(A_{[1,\infty)}^{k, k+1}\right)=0$, 
and, hence,  
$\Px\left(B_k\right)\geq \epsilon.$

\item[ 4)] If $r>z_2$ is strict, then, with positive probability, the maximal rate relocates as follows.
There exists $ \varepsilon> 0$ such that
\begin{equation}
\label{l51bound}
 \Px \left(B_k,    \max\limits_{i=k+2, k+3}\Gamma_{i}(X(\tau_{k+2}))
   = \max\limits_{i}\Gamma_{i}(X(\tau_{k+2})) \right) \geq \varepsilon,
\end{equation} 
where $\max_{i}\Gamma_i(X(\tau_{k+2}))$ may be attained at $k+3$ only if $\lambda_{k+3} > \lambda$.
\end{enumerate}
\end{lem}

 Part 4) of Lemma~\ref{L5} is similar to Lemmas~\ref{L2}-\ref{L3} in that it also
describes relocation of the maximal rate to a site with larger parameter $\lambda$. The main 
difference is that in Lemma~\ref{L5} the relocation time is 
random. This is in contrast to  Lemmas~\ref{L2}-\ref{L3}, where the relocation time is
 deterministic. 

The proposition and definition below are intended to clarify some assumptions and simplify some notations in Lemmas \ref{L6}-\ref{L8} below.

\begin{prop}
\label{P0}
Let $\{k, k+1\}$ be a local minimum of size 2 with $\lambda=\lambda_{k}=\lambda_{k+1}$, and let $r=r(\bx)$, $z_1$ and $z_2$ be quantities as in Definition (\ref{D4}).
Then, $z_1<z_2$ if and only if local minimum $\{k, k+1\}$ is of type 1, in which case there might  exist  $\bx$ such that $z_1<r<z_2$.
 Otherwise, if a local minimum $\{k, k+1\}$ is of type 2, then  $z_2\leq z_1$, in which case $r \geq z_2$ or $r \leq  z_1$ for all $\bx$.

\end{prop}

\begin{dfn}
\label{D5}
Remind that $\tau_k:=\inf(n: X_k(n)=X_k(0)+1)$ and let us further define the following stopping times
\begin{align*}
\sigma_k&=\min(\tau_{k-1}, \tau_{k+2}),\\
w_k&=\min(\tau_i \, : \, i\neq k\pm1, k, k+2),
\end{align*}
and following events  
\begin{align*}
D_k&=\{\sigma_k<w_k\},\\
D_k'&=\{\tau_{k-1}<\min(\tau_{k+2}, w_k)\},\\
D_k''&=\{\tau_{k+2}<\min(\tau_{k-1}, w_k)\}.
\end{align*}
\end{dfn}
Note that $D_k'\cap D_k''=\emptyset$, $D_k=D_k'\cup D_k''$ and  $A^{k, k+1}_{[1, \infty)}\cap D_k=\emptyset$. 

\begin{lem}
\label{L6}
Suppose that 
$
\{k, k+1\}$ is   a local minimum of size $2$,  and
$ \bx\in  \Z_{+}^N$ is such that
$\max(\Gamma_{k}(\bx), \Gamma_{k+1}(\bx))=\max_{i}\Gamma_i(\bx).$
\begin{enumerate}
\item[1)] There exists  $\epsilon>0$ such that  
$$
\Px\left(A_{[1,\infty)}^{k, k+1}\bigcup D_k\right)
=\Px\left(A_{[1,\infty)}^{k, k+1}\right)+\Px\left(D_k\right)\geq \epsilon.
$$
\item[2)]
If $ z_1<r< z_2$ (only possible if $\{k,k+1\}$ is of type 1),
then, with positive probability, all subsequent particles are allocated at sites $\{k, k+1\}$, i.e.  
$\Px\left(A_{[1,\infty)}^{k, k+1}\right)>\varepsilon$ for some  $\varepsilon>0$.

\item[3)]  If $r\leq z_1$ or $r\geq z_2$ (always the case if $\{k,k+1\}$ is of type 2),\\
then $\Px\left(A_{[1,\infty)}^{k, k+1}\right)=0$ and, hence, $\Px\left(D_k\right)\geq \epsilon$.

\end{enumerate}
\end{lem}

Lemma~\ref{L6}  is analogous to Parts 1)-3) of Lemma~\ref{L5} for the case of a local minimum of size 2.
An analogue of Part 4) of Lemma~\ref{L5} in the  same situation   is provided by Lemma~\ref{L7} below. 

\begin{lem}
\label{L7}
Suppose that local minimum  $\{k,k+1\}$ is of size $2$  with
$\lambda_k=\lambda_{k+1}:=\lambda$, and $ \bx\in  \Z_{+}^N$ is such that
$\max(\Gamma_{k}(\bx), \Gamma_{k+1}(\bx))=\max_{i}\Gamma_i(\bx).$

\begin{itemize}
\item[ 1)]  
If $\{k,k+1\}$ is of type $1$ and $r<z_1$, or $\{k,k+1\}$ is of type $2$ and $r<z_2$ then
\[
  \Px\left( D_k',
  \max\limits_{i=k-2,k-1}\Gamma_{i}(X(\tau_{k-1}))
 =\max\limits_{i=1,\ldots,N}\Gamma_{i}(X(\tau_{k-1}))\right) \geq \varepsilon>0.
\]
where $\max\limits_{i}\Gamma_i(X(\tau_{k-1}))$ may be attained at $k-2$ only if $\lambda_{k-2} > \lambda$.
\item[ 2)] 
If $\{k,k+1\}$ is of type $1$ and $r>z_2$, or $\{k,k+1\}$ is of type $2$ and $r>z_1$ then
\[
  \Px\left(D_k'',
  \max\limits_{i=k+2, k+3}\Gamma_{i}(X(\tau_{k+2}))
 =\max\limits_{i=1,\ldots,N}\Gamma_{i}(X(\tau_{k+2}))\right) \geq \varepsilon>0,
\]
where $\max\limits_{i}\Gamma_i(X(\tau_{k+2}))$ may be attained at $k+3$ only if $\lambda_{k+3} > \lambda$.
\item[ 3)] 
 If $\{k,k+1\}$ is of type $2$ and $z_2<r<z_1$, then  
\begin{multline*}
  \Px\left(D_k',
  \max\limits_{i=k-2, k-1}\Gamma_{i}(X(\tau_{k-1}))
 =\max\limits_{i=1,\ldots,N}\Gamma_{i}(X(\tau_{k-1}))\right) \\
+
 \Px\left(D_k'',
  \max\limits_{i=k+2, k+3}\Gamma_{i}(X(\tau_{k+2}))
 =\max\limits_{i=1,\ldots,N}\Gamma_{i}(X(\tau_{k+2}))\right)
\geq \varepsilon>0,
\end{multline*}
where $\max \Gamma_i$ follows the corresponding prescriptions as above.
\end{itemize}
\end{lem}

\begin{rmk}
\label{R1}
{\rm The next lemma concerns the borderline cases in between having a local minimum $\{k, k+1\}$ of size 2 and type 1 or a saddle point.
For example, in notations of Lemma~\ref{L7} 
these cases are formally obtained by setting either $\lambda_{k-1}=\lambda$ (where $-\infty=z_1<z_2$), or $\lambda_{k+2}=\lambda$ (where  $z_1<z_2=\infty$). 
As both cases can be addressed in similar ways, the lemma below deals only with the case $\lambda_{k-1}=\lambda$.}
\end{rmk}

\begin{lem}
\label{L8}
Suppose that sites  $\{k-1, k, k+1, k+2\}$  are such that 
$$\lambda_{k-1}=\lambda_k=\lambda_{k+1}=:\lambda<\lambda_{k+2},$$
$\bx\in  \Z_{+}^N$ is such that 
$\max(\Gamma_{k}(\bx), \Gamma_{k+1}(\bx))=\max_{i}\Gamma_i(\bx) $
  and, additionally,  
$\Gamma_{k-1}(\bx)\leq \Gamma_{k+1}(\bx)$.
\begin{enumerate}

\item[1)]  There exists  $\epsilon>0$ such that  
$$
\Px\left(A_{[1,\infty)}^{k, k+1}\bigcup D_k\right)
=\Px\left(A_{[1,\infty)}^{k, k+1}\right)+\Px\left(D_k\right)\geq \epsilon.
$$

\item[ 2)] 
 If $r<z_2$, 
then,  with positive probability all subsequent particles are allocated at sites $\{k, k+1\}$, i.e.  
$\Px\left(A_{[1,\infty)}^{k, k+1}\right)\geq \epsilon$  for some  $\varepsilon>0$.

\item[ 3)] 
If  
$r\geq z_2$,
then  $\Px\left(A_{[1,\infty)}^{k, k+1}\right)=0$ and, hence, $\Px\left(D_k\right)\geq \epsilon$.

\item[ 4)] If $r>z_2$, then there exists $ \varepsilon> 0$ such that

\begin{equation*}
 \Px \left(B_k,    \max\limits_{i=k+2, k+3}\Gamma_{i}(X(\tau_{k+2}))
   = \max\limits_{i}\Gamma_{i}(X(\tau_{k+2})) \right) \geq \varepsilon,
\end{equation*} 
where $\max\limits_{i}\Gamma_i(X(\tau_{k+2}))$ may be attained at $k+3$ only if $\lambda_{k+3} > \lambda$.

\end{enumerate}
\end{lem}

The following corollary concerns  those cases   covered by Parts 3) of Lemmas~\ref{L5}, \ref{L6} and \ref{L8},
where  the configuration parameter $r$ is  equal to one of the model parameters $z_1$ and $z_2$. In what follows we call 
them critical cases.

\begin{cor}
\label{critical}
For the critical cases, relocation of the maximal rate to a site with larger parameter $\lambda$ also occurs, with positive probability, in finite time. 
\end{cor}

\begin{rmk}
\label{R2}
{\rm Let us  remark the following.

\begin{itemize}
\item[1)]
It is important to emphasize that in all the above cases where the maximal rate $\max_i\Gamma_i(\bx)$ eventually relocates with positive probability, it always relocates to a site with strictly larger parameter  $\lambda$.
\item[ 2)] Note that  Lemmas~\ref{L2}, \ref{L3}, \ref{L5} and \ref{L7} can be appropriately reformulated in order to cover the symmetric cases by simply re-labelling the graph sites in reverse order as the graph is a cycle. 
For example, if  $\{k, k+1\}$ is a saddle point as in Lemma~\ref{L5}, then the corresponding symmetric case would be $\lambda_{k-1}>\lambda_k=\lambda_{k+1}>\lambda_{k+2}$, etc.
\end{itemize}
 }
\end{rmk}

\section{Random geometric progressions and\\
 Bernoulli measures}
\label{App} 

The statements and propositions in this section are essential building blocks for the proof of lemmas which follow.
The reason is that along the proofs of Lemmas~\ref{L4}-\ref{L8} we need to analyse the limiting behaviour of random variables of the form $\sum_{i=0}^{n} \prod_{j=1}^i\zeta_j$,
as $n \rightarrow \infty$,  where $\{\zeta_j, \> j \geq 1\}$ is an i.i.d. sequence of positive random variables. It will also be necessary to compare such variables and introduce some stochastic domination concepts to enable us to carry out uniform estimates not depending on the starting configuration $X(0) = \bx$.
We refer to \cite{Tho2000} for standard definitions and basic properties of stochastic domination.
The following notations are used throughout. Given random variables $X$ and $Y$ (or sequences $X$ and $Y$), we write $X\geq_{st}Y$ if $X$ stochastically dominates $Y$.
Similarly, given two probability measures $\nu$ and $\mu$, we write $\mu\geq_{st}\nu$ if $\mu$ stochastically dominates $\nu$.

\paragraph{Random geometric progressions.}
In this subsection we consider random variables realised on a certain probability space $(\Omega, {\cal F}, \P)$.
$\E$ denotes the expectation with respect to probability measure $\P$.
If $X$ and $Y$ are random variables or sequences such that $X\geq_{st}Y$, then we may assume that  $\P$ is a coupling of probability distributions of $X$ and $Y$ such that  $\P(X\geq Y)=1$. 
Such a coupling exists by Strassen's theorem (\cite{Str65}).

Given   a random    sequence  ${\bf \zeta}=\{\zeta_i,\, i\geq1\}$, define 
\begin{equation}
\label{sequences}
Y_i=\prod_{j=1}^i\zeta_j, i\geq 1,\, Y_0=1,\,\, \mbox{and}\,\, Z_n(\zeta)=\sum_{i=0}^{n}Y_i,\, n\geq 1,
\end{equation}
and
\begin{equation}
\label{Z_infty}
Z({\bf \zeta})=\sum_{i=0}^{\infty}Y_i.
\end{equation}
\begin{prop}
\label{P2}
1) 
Let ${\bf \zeta}=\{\zeta_{ i},\, i\geq1\}$ be an i.i.d. sequence of positive random variables such that
$\E\left(\log(\zeta_{i})\right)<0$. 
Then 
$\P(Z(\zeta)<\infty)=1$ 
and, consequently, 
$\E\left(e^{-Z(\zeta)}\right)>0$.

2) Let ${\bf \theta}=\{\theta_{i},\, i\geq 1\}$ be  another i.i.d. sequence of positive random variables such 
that  $\E\left(\log(\theta_{i})\right)<0$ and  ${\bf \theta}\geq_{st}\zeta$. Then 
$\E\left(e^{-Z(\zeta)}\right)\geq \E\left(e^{-Z(\theta)}\right).$
\end{prop}

{\it Proof of Proposition~\ref{P2}.}
Denote $\E\left(\log(\zeta_{i})\right)=a<0$. Given $\delta>0$ such that $a+\delta<0$, it follows from the strong law of large numbers that 
$Y_n<e^{(a+\delta)n}$
for all but finitely many $n$ almost surely. Therefore, a tail of $Z(\zeta)$ is eventually  majorised by 
the corresponding tail of a converging geometric progression. 
In turn, finiteness of $Z(\zeta)$ implies positiveness of the expectation.
Moreover, note that $e^{Z(\cdot)}$ is an increasing function.
Therefore, $e^{-Z(\zeta)}\geq_{st}e^{-Z(\theta)}$ and hence, 
$\E\left(e^{-Z(\zeta)}\right)\geq \E\left(e^{-Z(\theta)}\right)$ as claimed.
\qed

\begin{dfn}
\label{reciprocal}
Let $\zeta=\{\zeta_i,\, i\geq1\}$ and $\eta=\{\eta_j, j\geq 1\}$ be i.i.d. sequences of positive random variables. 
Sequence $\eta$ is said to be reciprocal to $\zeta$ if $\eta_1$ has the same distribution as  $1/\zeta_{1}$. 
\end{dfn}

The following proposition  follows  from  basic properties of stochastic domination.
\begin{prop}
\label{dominance}
Let 
 $X$ and  $Y$ be  two i.i.d.  sequences of  positive random variables,  and let  $\eta_{\scriptscriptstyle X}$ and $\eta_{\scriptscriptstyle Y}$
  be their corresponding reciprocal sequences. If $X\geq_{st}Y$ then 
 $\eta_{\scriptscriptstyle X}\leq_{st}\eta_{\scriptscriptstyle Y}$.
\end{prop}

\begin{prop}
\label{F}
Let $\zeta=\{\zeta_{ i},\, i\geq1\}$ be an i.i.d. sequence of  positive random variables such that 
$\E\left(\log(\zeta_{i})\right)>0$. 
Let  $\{Y_i, i\geq 0\} $ and $\{Z_{n}(\zeta), n\geq 1\}$ be the random variables as in (\ref{sequences}).
Define the following random sequence
\begin{equation*}
F_{n}(\zeta)=Z_{n}(\zeta)/Y_{n},\, n\geq 1.
\end{equation*}
Then, $F_{n}(\zeta)$ converges in distribution to 
\begin{equation}
\label{F=Z}
Z(\eta)=1+\sum\limits_{i=1}^{\infty}\prod_{j=1}^i\eta_j, \quad \text{as} \>\>\> n\to \infty,
\end{equation}
where $\eta$ is the sequence reciprocal to $\zeta$. Moreover, $Z(\eta)$ is almost surely finite and $Z(\eta)\geq_{st}F_n(\zeta)$ for any $n\geq 1$.
\end{prop}
{\it Proof of Proposition~\ref{F}.}
First, note that for every $n\geq 1$,
$$F_{n}(\zeta)=1+\sum\limits_{i=1}^{n}\prod_{j=1}^i\zeta_{n-j+1}^{-1}= 1+\sum\limits_{i=1}^{n}\prod_{j=1}^i\eta_{j}^{(n)},$$
where $\eta_{j}^{(n)}=\zeta_{n-j+1}^{-1}$.
This means that  $F_{n}(\zeta)$ has the same distribution as $Z_{n}(\eta)$ defined for the sequence $\eta=\{\eta_i,\, i\geq 1\}$ reciprocal to $\zeta$.
Therefore, $F_n(\zeta)$ converges in distribution to $Z(\eta)$.
In addition, $\E(\log(\eta_1))=-\E(\log(\zeta_1))<0$.
Therefore, by Proposition \ref{P2}, $Z(\eta)$ is almost surely finite. 
Finally, it follows by construction that  $Z(\eta)\geq_{st}F_n(\zeta)$, $n\geq 1$.
\qed

\begin{prop}
\label{P3}
Let $\zeta=\{\zeta_{ i},\, i\geq1\}$ be an i.i.d. sequence of positive random variables such that 
$\E\left(\log(\zeta_{i})\right)=a>0$, and $\eta=\{\eta_i,\, i\geq 1\}$ be its reciprocal sequence.
Given  $0<\gamma<1$, define the following stopping time
\begin{equation}
\label{hat_m}
\widehat m=\min(n: \gamma Y_n\geq 1).
\end{equation}
Then both  $Z(\eta)<\infty$ and $Z_{\widehat m-1}(\zeta)<\infty$  almost surely,   
$\gamma Z_{\widehat m-1}(\zeta) \leq_{st} Z(\eta)$, and, hence,   
\begin{equation}
\label{ex1}
\E\left(e^{-\gamma Z_{\widehat m-1}(\zeta)}\right)\geq \E\left(e^{-Z(\eta)}\right)>0.
\end{equation} 
 \end{prop}
{\it Proof of Proposition \ref{P3}.} 
By Proposition \ref{F}, $Z(\eta)$ is almost surely finite and $F_n(\zeta)\leq_{st} Z(\eta)$ for all $n\geq 1$.
Therefore, $F_{\widehat m-1}(\zeta)\leq_{st}Z(\eta)$.
Since $\gamma Y_{\widehat m-1}<1$ we obtain that $$\gamma Z_{\widehat m-1}(\zeta)<Z_{\widehat m-1}(\zeta)/Y_{\widehat m-1}= F_{\widehat m-1}(\zeta).$$
Consequently, $\gamma Z_{\widehat m-1}(\zeta) \leq_{st} Z(\eta)$, which implies (\ref{ex1}) as claimed. 
\qed

\begin{prop}
\label{C3}
Let $\zeta=(\zeta_{i},\, i\geq 1)$ and $\theta=(\theta_{i},\, i\geq 1)$ be i.i.d. sequences of positive random variables such that both $\E\left(\log(\zeta_{1})\right)>0$ and $\E\left(\log(\theta_{1})\right)>0$, and  $\zeta_1\geq_{st}\theta_1$. 
Let $\eta_{\scriptscriptstyle\zeta}$ and $\eta_{\theta}$ be sequences reciprocal to $\zeta$ and $\theta$, respectively.
Given $0<\gamma<1$, let $\widehat m$ be the stopping time for sequence $\zeta$ as in (\ref{hat_m}).
Then $$\E\left(e^{-\gamma Z_{\widehat m-1}(\zeta)}\right)\geq \E\left(e^{-Z(\eta_{\theta})}\right).$$
\end{prop}
{\it Proof of Proposition  \ref{C3}.}
By Proposition \ref{F} both $Z(\eta_{\zeta})$ and $Z(\eta_{\theta})$ are almost surely finite.
Further, by Proposition \ref{dominance}  $\eta_{\zeta}\leq_{st}\eta_{\theta}$.
Therefore $$\E\left(e^{-Z(\eta_{\zeta})}\right)\geq \E\left(e^{-Z(\eta_{\theta})}\right).$$
By Proposition \ref{P3}, it follows that
$$\E\left(e^{-\gamma Z_{\widehat m-1}(\zeta)}\right)\geq \E\left(e^{-Z(\eta_{\zeta})}\right)\geq 
\E\left(e^{-Z(\eta_{\theta})}\right)$$
as claimed. \qed

\paragraph{Bernoulli measures.} 
Now, we introduce a family of Bernoulli measures and some notations that will be used throughout proofs of Lemmas~\ref{L4}-\ref{L8}. 

Let $\xi=(\xi_i, i\geq 1)$ be a sequence of  independent Bernoulli random variables with success probability $p$.
Let $\mu_p$ be the distribution of $\xi$, that is, the product Bernoulli measure defined 
on the set of  infinite binary sequences, and denote by $\Ee_p$ the expectation with respect to the Bernoulli measure $\mu_p$. 

Define 
\begin{equation}
\label{Ui}
U_i=\xi_1+\cdots + \xi_i,\, i\geq 1,
\end{equation}
the binomial random variables corresponding to  a  Bernoulli sequence $\xi$. 

Let $\lambda_{k-1}$, $\lambda_k$, $\lambda_{k+1}$ and  $\lambda_{k+2}$ be $\lambda$-parameters corresponding to quadruples $\{k-1, k, k+1, k+2\}$ of the graph sites such that $\lambda=\lambda_k=\lambda_{k+1}$ as in  Lemmas~\ref{L4}-\ref{L8}.
Let us define the following i.i.d. sequences
\begin{equation}
\label{zetas}
\begin{split}
\zeta_1&=(\zeta_{1, i}=e^{\lambda_{k-1}(1-\xi_i)-\lambda},\, i\geq 1),\\
\zeta_2&=(\zeta_{2, i}=e^{\lambda_{k+2}\xi_i-\lambda},\, i\geq 1).
\end{split}
\end{equation}
It is a  well known fact  that if $ 0<p'\leq p''<1$, then $\mu_{p'}\leq_{st}\mu_{p''}$. 
This fact yields the following proposition. 
\begin{prop}
\label{zeta(p1)<zeta(p2)}
Let $\zeta_1', \zeta_2'$ and $\zeta_1'', \zeta_2''$ be sequences defined by (\ref{zetas}) for Bernoulli sequences with success probabilities $p'$ and $p''$, respectively.
If $0<p'\leq p''<1$, then $\zeta_1'\geq_{st}\zeta_1''$ and $\zeta_2'\leq_{st}\zeta''_2$.
\end{prop}
Note that variables $Z_n$ (defined in  (\ref{sequences})) corresponding to sequences $\zeta_1$ and $\zeta_2$ can be 
 expressed in terms of Binomial random variables (\ref{Ui}) as follows 
\begin{equation}
\label{Z_for_zeta}
Z_n(\zeta_1)=\sum\limits_{i=0}^{n}e^{\lambda_{k-1}(i-U_i)-\lambda i}\quad\text{and}\quad
Z_n(\zeta_2)=\sum\limits_{i=0}^{n}e^{\lambda_{k+2}U_i-\lambda i}.
\end{equation}
It is useful to note  that if $\lambda_{k-1}=\lambda_{k+2}=\lambda$, then  the above expressions are
\begin{equation*}
Z_n(\zeta_1)=\sum\limits_{i=0}^{n}e^{-\lambda U_i} \quad\text{and}\quad
Z_n(\zeta_2)=\sum\limits_{i=0}^{n}e^{\lambda(U_i-i)}.
\end{equation*}

\section{Proofs of Lemmas}
\label{ProofsLemmas}

In the following proofs we show the existence of positive real constants $C$, $c$, $\epsilon$ and $\eps$, whose exact values are immaterial and may vary from line to line, but which do not depend on the starting configuration $X(0) = \bx $. In order to avoid notational clutter we shall denote $\Gamma_i(\bx)$ simply by $\Gamma_i$ for all $i$. Moreover, whenever we fix index $k \in \{1,\ldots,N\}$ and consider indices in the neighbourhood of $k$, those indices should be interpreted as $modulo \> N$.

\subsection{Proofs of Lemmas~\ref{L1}-\ref{L3}}
\label{prL1-L3}

For short, denote $B=\sum_{i \neq k, k \pm 1}\Gamma_i$ and $Z=\sum_{i=1}^N\Gamma_i$.
By assumption, $\Gamma_k=\max\limits_{i=1, \ldots, N}\Gamma_i$, then
\begin{equation}
\label{b1}
\frac{\Gamma_{k-1}}{\Gamma_k}\leq 1, \frac{\Gamma_{k+1}}{\Gamma_k}\leq 1, \, \Gamma_k\geq \frac{Z}{N},\, \text{and}\, \frac{Z-\Gamma_k}{Z}\leq \frac{(N-1)}{N}.
\end{equation}
It follows from the last two inequalities that 
\begin{equation}
\label{b2}
\frac{B}{\Gamma_k}\leq N-1.
\end{equation}

\paragraph{{\it Proof of Lemma~\ref{L1}.}}
Remind that $\lambda_k>\max(\lambda_{k-1}, \lambda_{k+1})$. We need to prove the existence of a positive number $\epsilon$ such that 
\begin{equation}
  \label{locmaxL1}
   \Probx\left(A_{[1, \infty)}^k\right) =
   \prod\limits_{n=0}^{\infty}
   \frac{\Gamma_ke^{\lambda_kn}}
   {\Gamma_{k-1}e^{\lambda_{k-1}n}+\Gamma_ke^{\lambda_kn}+\Gamma_{k+1}e^{\lambda_{k+1}n}+ B}
   >\epsilon,
\end{equation}
where $\epsilon> 0$ depends only on $\lambda_{k-1}, \lambda_k, \lambda_{k+1}$ and  $N$.

Indeed, rewriting the identity in (\ref{locmaxL1}) and applying bounds (\ref{b1}) and (\ref{b2}),  
\begin{align*}
  &\Probx\left(A_{[1, \infty)}^k\right)\\
  &=\exp\left(-\sum_{n=0}^{\infty}\log\left(1+\frac{\Gamma_{k-1}}{\Gamma_k}e^{(\lambda_{k-1}-\lambda_k)n}+
    \frac{\Gamma_{k+1}}{\Gamma_k}e^{(\lambda_{k+1}-\lambda_k)n}+\frac{B}{\Gamma_k}e^{-\lambda_kn}\right)\right) \\
  & \geq \exp\left(-\sum_{n=0}^{\infty}\log(1+e^{(\lambda_{k-1}-\lambda_k)n}+e^{(\lambda_{k+1}-\lambda_k)n}+(N-1)e^{-\lambda_kn})\right) \\
  & \geq     \exp\left(-C\sum_{n=0}^{\infty}(e^{(\lambda_{k-1}-\lambda_k)n}+e^{(\lambda_{k+1}-\lambda_k)n}+(N-1)e^{-\lambda_kn})\right)>\epsilon>0,
\end{align*}
since the series in the exponent above converges. It is not hard to see that in the last inequality, $\epsilon$ should depend only on $\lambda_{k-1}, \lambda_k, \lambda_{k+1}$ and $N$. \qed

\paragraph{{\it Proof of Lemma~\ref{L2}}.}
Remind that $\lambda_{k-1}<\lambda_k<\lambda_{k+1}$. We need to prove the existence of a finite positive integer $\hat n$ and a positive number $\epsilon$ such that 
\begin{equation}
\label{n1}
\Gamma_{k+1}e^{\lambda_{k+1}\hat n} \geq \Gamma_ke^{\lambda_k\hat n} > \max\left(\Gamma_{k-1}e^{\lambda_{k-1}\hat n}, \max_{i \neq k, k \pm 1}\Gamma_i\right)
\end{equation} 
and 
\begin{equation}
  \label{pntgrowthL2}
  \Probx(A_{[1, \hat n]}^k) = 
  \prod\limits_{n=0}^{\hat n}
  \frac{\Gamma_ke^{\lambda_kn}}
  {\Gamma_{k-1}e^{\lambda_{k-1}n}+\Gamma_ke^{\lambda_kn}+\Gamma_{k+1}e^{\lambda_{k+1}n}+ B}
  >\epsilon,
\end{equation}
  where $\epsilon > 0$ depends only on $\lambda_{k-1}, \lambda_k, \lambda_{k+1}$ and $N$.
Note  that the sequence $e^{(\lambda_{k+1}-\lambda_k)n}, \> n \geq 0$ is strictly increasing, 
so there exists the minimal integer $\hat n$ such that 
\[e^{(\lambda_{k+1}-\lambda_k)\hat n} \geq \frac{\Gamma_k}{\Gamma_{k+1}}, \quad \text{that is,} \quad \frac{\Gamma_{k+1}(\bx + \hat n \be_k)}{\Gamma_k(\bx + \hat n \be_k)} \geq 1.\]
Then, there exists a positive constant such that $(\Gamma_{k+1}/\Gamma_k)e^{(\lambda_{k+1}-\lambda_k)\hat n} \leq C_1,$
and, hence, 
\begin{equation}
\label{geom}
\frac{\Gamma_{k+1}}{\Gamma_k}\sum\limits_{n=0}^{\hat n} e^{(\lambda_{k+1}-\lambda_k)n}\leq C_2<\infty,
\end{equation}
where $C_2$ depends only on $\lambda_k$ and $\lambda_{k+1}$.
Further, rewriting the identity in (\ref{pntgrowthL2}) and using bounds (\ref{b1}), (\ref{b2})  and (\ref{geom}), gives that

\begin{align*}
&\Probx(A_{[1, \hat n]}^k) \\
&=\exp\left(-\sum_{n=0}^{\hat n}\log\left(1+\frac{\Gamma_{k-1}}{\Gamma_k}e^{(\lambda_{k-1}-\lambda_k)n} + \frac{\Gamma_{k+1}}{\Gamma_k}e^{(\lambda_{k+1}-\lambda_k)n}+\frac{B}{\Gamma_2}e^{-\lambda_k n}\right)\right) \\
&\geq \exp\left(-\sum_{n=0}^{\hat n}\log\left(1+e^{(\lambda_{k-1}-\lambda_k)n}+\frac{\Gamma_{k+1}}{\Gamma_k}e^{(\lambda_{k+1}-\lambda_k)n}+(N-1)e^{-\lambda_kn}\right)\right) \\
& \geq \exp\left(-C_3\sum_{n=0}^{\hat n}\left(e^{(\lambda_{k-1}-\lambda_k)n}+\frac{\Gamma_{k+1}}{\Gamma_k}e^{(\lambda_{k+1}-\lambda_k)n}+(N-1)e^{-\lambda_k n}\right)\right)>\epsilon>0.
\end{align*}
\qed

\paragraph{\it Proof of Lemma~\ref{L3}.}
Recall  that $\lambda_k < \min(\lambda_{k-1}, \lambda_{k+1})$. As in the proof of Lemma~\ref{L2}, 
we need to  show  existence of a finite positive integer $\hat n$ and a positive $\epsilon$ such that 
\begin{equation}
\label{n2}
\max(\Gamma_{k-1}e^{\lambda_{k-1}\hat n}, \Gamma_{k+1}e^{\lambda_{k+1}\hat n}) \geq \Gamma_ke^{\lambda_k\hat n} \geq \max_{i \neq k, k \pm 1}\Gamma_i
\end{equation}
 and  
 \[
  \Probx(A_{[1, \hat n]}^k) =  
  \prod\limits_{n=0}^{\hat n}
  \frac{\Gamma_ke^{\lambda_kn}}
  {\Gamma_{k-1}e^{\lambda_{k-1}n}+\Gamma_ke^{\lambda_kn}+\Gamma_{k+1}e^{\lambda_{k+1}n}+ B}
  >\epsilon,
\]
where $\epsilon > 0$ depends only on $\lambda_{k-1}, \lambda_k, \lambda_{k+1}$ and  $N$.
This can be shown similar to  the proof of Lemma~\ref{L2}, and we skip details.  \qed

\subsection{Proofs of Lemmas~\ref{L4}-\ref{L8}}
\label{prL5-L8}

\subsubsection{Notations}
We start with some  preliminary considerations  and notations  that will be used throughout the proofs of Lemmas~\ref{L4}-\ref{L8}.

Let $\{k, k+1\}$  be a pair of sites such  that $\lambda_k=\lambda_{k+1}=\lambda$.
If, as defined in Definition~\ref{D2}, $r=r(\bx)=x_{k+2}-x_{k-1}$, then $\frac{\Gamma_{k+1}(\bx)}{\Gamma_{k}(\bx)}=e^{\lambda r}$. 
Therefore,  given that the next particle is allocated at either  $k$ or $k+1$, the conditional $\P_{\bx}$-probability to choose $k+1$ is equal to 
\begin{equation}
\label{p}
p:=p(r)=\frac{\Gamma_{k+1}(\bx)}{\Gamma_{k}(\bx)+\Gamma_{k+1}(\bx)}=\frac{e^{\lambda r}}{1+e^{\lambda r}}.
\end{equation}
We henceforth denote  $q=1-p$.
Furthermore, probability $p$ does not change by adding particles at sites $k$ and $k+1$ since configuration parameter $r$ remains constant. 

Note that $p(z)$, considered as a function of $z\in \R$, is monotonically increasing. 
A direct computation gives that unique solutions of equations 
$\lambda_{k-1}-\lambda=p(z)\lambda_{k-1}$ and $\lambda_{k+2}p(z)=\lambda$ 
are quantities $z_1$ and $z_2$ (defined in (\ref{z})) respectively. 

Let $S_n$ be the number of particles at site $k+1$ at time $n \geq 1$.
Let $S_0=0$ and $s(n)=(s_0, s_1, \ldots, s_n)$ be a fixed  trajectory of a finite  random sequence $S(n)=(S_0, S_1, \ldots, S_n)$.
Note that, by construction, any trajectory $s(n)$ is a sequence of non-negative integers  such that $s_0=0$ and $s_{i}-s_{i-1}\in\{0, 1\}$, $i=1, \ldots, n$.

For short, denote 
\begin{align*}
\nonumber
\Gamma_i&=\Gamma_i(\bx),\, \widetilde \Gamma_k=\sum_{i\neq k, k\pm1, k+2}\Gamma_i,\\
\gamma_{k,1}&=\frac{\Gamma_{k-1}}{\Gamma_{k}+\Gamma_{k+1}},\, 
\gamma_{k,2}=\frac{\Gamma_{k+2}}{\Gamma_{k}+\Gamma_{k+1}}, \, \widetilde \gamma_k=\frac{\widetilde \Gamma_k}{\Gamma_{k}+\Gamma_{k+1}}.
\end{align*}

In the rest of this section we are going to derive expressions for probabilities $\Px\left(A^{k, k+1}_{[1, n+1]}\right), \> n \geq 1$, in terms of expectations  with respect to a Bernoulli product measure on $\{0,1\}^{\infty}$ with parameter $p$ defined in (\ref{p}). 
These expressions allow one to obtain lower and upper bounds for the above probabilities.
We start with the case of fixed $n$ and then extend it to the case where $n$ is a stopping time.

In the above  notations
 \begin{align*}
\Px\left(A^{k, k+1}_{i+1}, S_{i+1}\bigg.
\right.&\left.=s_{i+1}\bigg|A^{k, k+1}_{[1, i]}, S_i=s_i\right)\\
&=
\frac{p^{s_{i+1}-s_i}q^{1-(s_{i+1}-s_i)}(\Gamma_{k}+\Gamma_{k+1})e^{\lambda i}}{(\Gamma_{k}+\Gamma_{k+1})e^{\lambda i}
+\Gamma_{k-1}e^{\lambda_{k-1}(i-s_i)}+\Gamma_{k+2}e^{\lambda_{k+2}s_i}+\widetilde \Gamma_k}\\
&=\frac{p^{s_{i+1}-s_i}q^{1-(s_{i+1}-s_i)}}{1+\gamma_{k,1}e^{\lambda_{k-1}(i-s_i)-\lambda i}+\gamma_{k,2}e^{\lambda_{k+2}s_i-\lambda i}+
\widetilde \gamma_ke^{-\lambda i}}.
\end{align*}
Then,  given $n$ we obtain by repeated conditioning that  
\begin{multline*}
\Px\left(A^{k, k+1}_{[1, n+1]},  S_{n+1}=s_{n+1}, \ldots, S_1=s_1\right)\\
=p^{s_{n+1}-s_n}q^{1-(s_{n+1}-s_n)}p^{s_n}q^{n-s_n}W_n(s_1, \ldots, s_n),
\end{multline*}
where 
\begin{equation}
\label{W}
W_n(s_1, \ldots, s_n)=
\prod\limits_{i=0}^n \frac{1}{1+\gamma_{k,1}e^{\lambda_{k-1}(i-s_i)-\lambda i}+
\gamma_{k,2}e^{\lambda_{k+2}s_i-\lambda i}+\widetilde\gamma_ke^{-\lambda i}}.
\end{equation}
Consequently, we get that  
\begin{equation}
\label{total}
\begin{split}
\Px\left(A^{k, k+1}_{[1, n+1]}\right)&=\sum\limits_{s(n+1)}p^{s_{n+1}}q^{n+1-s_{n+1}}W_n(s_1, \ldots, s_n),\\
&=\sum\limits_{s(n)}(p+q)p^{s_{n}}q^{n-s_{n}}W_n(s_1, \ldots, s_n),\\
&=\sum\limits_{s(n)}p^{s_{n}}q^{n-s_{n}}W_n(s_1, \ldots, s_n),
\end{split}
\end{equation}
where the  sum in the first line is over all possible trajectories 
$s(n+1)=(s_1, \ldots, s_{n+1})$ of $S(n+1)=(S_1, \ldots, S_{n+1})$ and the other two 
are over all possible trajectories 
$s(n)=(s_1, \ldots, s_n)$ of $S(n)=(S_1, \ldots, S_n)$.
Therefore, we arrive to the following 
equation
\begin{equation}
\label{random}
\Px\left(A^{k, k+1}_{[1, n+1]}\right)=\Ee_p(W_n(U_1, \ldots, U_n)),
\end{equation}
where
  $\Ee_p$ is the expectation  with respect to the  Bernoulli measure $\mu_p$ defined in Section \ref{App} and 
$U_i$, $i\geq1$, are Binomial random variables defined in (\ref{Ui}).

Further, assumptions of  Lemmas \ref{L4}-\ref{L8} imply that 
$\frac{\Gamma_i}{\Gamma_k+\Gamma_{k+1}}\leq 1$, $i=1, \ldots, N$.
Therefore, 
\begin{equation}
\label{reminder}
\widetilde \gamma_ke^{-\lambda i}\leq (N-4)e^{-\lambda i}\leq c_1e^{-c_2i},
\end{equation}
for some $c_1, c_2>0$,  as $s_i\leq i$, $i=1, \ldots, n$.
Using bound  (\ref{reminder}) and  inequality  $\log(1+z)\leq z$ for all $z\geq 0$ we obtain that 
\begin{align}
\nonumber
W_n(s_1, \ldots, s_n)&\geq \prod_{i=0}^n
\frac{1}{1+
\gamma_{k,1}e^{\lambda_{k-1}(i-s_i)-\lambda i}+
\gamma_{k,2}e^{\lambda_{k+2}s_i-\lambda i}+c_1e^{-c_2i}}\\
&=e^{-\sum_{i=0}^n\log\left(1+\gamma_{k,1}e^{\lambda_{k-1}(i-s_i)-\lambda i}+
\gamma_{k,2}e^{\lambda_{k+2}s_i-\lambda i}+c_1e^{-c_2i}\right)}\nonumber\\
&\geq e^{-\left(\sum_{i=0}^n \gamma_{k,1}e^{\lambda_{k-1}(i-s_i)-\lambda i}+
\gamma_{k,2}e^{\lambda_{k+2}s_i-\lambda i}+c_1e^{-c_2i}\right) }\nonumber\\
&\geq \delta e^{-\gamma_{k,1}\sum_{i=0}^n   e^{\lambda_{k-1}(i-s_i)-\lambda i}}
e^{-\gamma_{k,2}\sum_{i=0}^n  e^{\lambda_{k+2}s_i-\lambda i}}, \label{27}
\end{align}
for some $\delta>0$ not depending on the configuration $\bx$.
On the other hand, note that   
\begin{equation}
\label{W1}
W_n(s_1, \ldots, s_n)\leq 
\prod\limits_{i=0}^n \frac{1}{1+
\gamma_{k,1}e^{\lambda_{k-1}(i-s_i)-\lambda i}+
\gamma_{k,2}e^{\lambda_{k+2}s_i-\lambda i}}.
\end{equation}
The above inequalities yield the following lower and upper  bounds  
\begin{equation}
\label{lower}
\Px\left(A^{k, k+1}_{[1, n+1]}\right)\geq 
\delta \Ee_p\left( e^{-\gamma_{k,1}\sum_{i=0}^n   e^{\lambda_{k-1}(i-U_i)-\lambda i}}
e^{-\gamma_{k,2}\sum_{i=0}^n  e^{\lambda_{k+2}U_i-\lambda i}}  \right),
\end{equation}
\begin{equation}
\label{upper}
\Px\left(A^{k, k+1}_{[1, n+1]}\right)
\leq \Ee_p\left(
\prod\limits_{i=0}^n \frac{1}{1+
\gamma_{k,1}e^{\lambda_{k-1}(i-U_i)-\lambda i}+
\gamma_{k,2}e^{\lambda_{k+2}U_i-\lambda i}}\right).
\end{equation}
We will also need a
 generalisation of  lower bound (\ref{lower})
for probabilities  $\Px\left(A^{k, k+1}_{[1, \tau]}\right)$, where 
$\tau$ is one of the following stopping times, 
$\min(n: S_n-c_1n\geq c_2)$,  $\min(n: n-S_n\geq c_3)$,
and the minimum of two such stopping times. 
At the moment, we shall not further specify such stopping times as
it  will be clear later which one it refers to. 
Arguing similarly as in equation (\ref{total}), one can obtain that  
\begin{equation}
\label{total0}
\Px\left(A^{k, k+1}_{[1, \tau]}\right)=
\sum\limits_{n=0}^{\infty}
\sum\limits_{s(n)}p^{s_n}q^{n-s_n}W_n(s_1, \ldots, s_n)1_{\{\mathsf M_n\}},
\end{equation}
where  $\mathsf M_n$ is a set of paths $s(n)=(s_1, \ldots, s_n)$ for which $\tau=n+1$.
Furthermore, similar to equation (\ref{random}), we can rewrite equation above as 
\begin{equation}
\label{random2}
\Px\left(A^{k, k+1}_{[1, \tau]}\right)=
\Ee_p\left(W_{\tilde \tau}(U_1, \ldots, U_{\tilde \tau})\right),
\end{equation} 
where
$\tilde \tau$ is a stopping time defined by replacing $S_n$ by $U_n$ in the same way as $\tau$ but in terms of random variables $U_n$.
Proceeding similar to how we got lower bound (\ref{lower}) we obtain  the following lower bound
\begin{equation}
\label{Bern31}
\Px\left(A^{k, k+1}_{[1, \tau]}\right)
\geq \delta \Ee_p\left( e^{-\gamma_{k,1}\sum_{i=0}^{\tilde \tau-1} e^{\lambda_{k-1}(i-U_i)-\lambda i}}
e^{-\gamma_{k,2}\sum_{i=0}^{\tilde  \tau-1}  e^{\lambda_{k+2}U_i-\lambda i}} \right).
\end{equation}

Let us rewrite the lower bounds in terms of random sequences $\zeta_1$, $\zeta_2$ and $Z_n$ as defined in (\ref{zetas}) and  (\ref{Z_for_zeta}). 
In these  notations, lower bounds (\ref{lower}) and (\ref{Bern31}) take the following form 
\begin{equation}
\label{lower_Z}
\Px\left(A^{k, k+1}_{[1, n+1]}\right)
\geq \delta \Ee_p\left( e^{-\gamma_{k,1}Z_{n}(\zeta_1)}
e^{-\gamma_{k,2}Z_{n}(\zeta_2)} \right)
\end{equation}
and 
\begin{equation}
\label{Bern31_Z}
\Px\left(A^{k, k+1}_{[1, \tau]}\right)
\geq \delta \Ee_p\left( e^{-\gamma_{k,1}Z_{\tilde\tau-1}(\zeta_1)}
e^{-\gamma_{k,2}Z_{\tilde\tau-1}(\zeta_2)} \right)
\end{equation}
respectively. 

Finally,  letting $n \rightarrow \infty$ in (\ref{lower}) and (\ref{lower_Z}) we obtain the following bound  
\begin{equation}
\label{lower1}
\begin{split}
\Px\left(A^{k, k+1}_{[1, \infty)}\right)&\geq 
\delta \Ee_p\left( e^{-\gamma_{k,1}\sum_{i=0}^{\infty} e^{\lambda_{k-1}(i-U_i)-\lambda i}}
e^{-\gamma_{k,2}\sum_{i=0}^{\infty}e^{\lambda_{k+2}U_i-\lambda i}}  \right)\\
&=
\delta \Ee_p\left( e^{-\gamma_{k,1}Z(\zeta_1)}
e^{-\gamma_{k,2}Z(\zeta_2)} \right).
\end{split}
\end{equation}

\subsubsection{Proof of Lemma \ref{L4}}
\label{prL4}
We start with the following proposition.
\begin{prop}
\label{PL4}
Let  $\mu_p$ be the Bernoulli measure  defined in Section~\ref{App},  and let  $U_n,\, n\geq 1,$ be the corresponding 
Binomial random variables (defined in (\ref{Ui})). Then
\begin{itemize}
\item[1)] 
given  $\eps \in (0,1)$ and $\kappa > 0$, there exist positive constants $c_1 $ and $ c_2$ such that  
\begin{equation}
\label{PL4_1}
\inf_{p \in (0,1)} \mu_p\left(\bigcap_{n=M}^{\infty}\left\{\frac{n}{2}p(1-\kappa)-c_1\leq U_n\leq np(1+\kappa)+c_2\right\}\right)\geq \eps,
\end{equation}
where $M=[p^{-1}]$ is the integer part of $p^{-1}$;
\item[2)] 
given $\lambda>0$, there exists $\eps_1>0$ such that 
$$\inf\limits_{p\in (0, 1)}\Ee_p\left(e^{-p\sum_{i=0}^{\infty}e^{-\lambda U_i}}\right)\geq \eps_1.$$
\end{itemize}
\end{prop}
{\it Proof of Proposition \ref{PL4}.}
Set $U_0=0$ and define the following random variables 
\begin{align*}
V_j&=U_{jM}-U_{(j-1)M}=\sum_{i=(j-1)M+1}^{jM}\xi_{i},\, j\geq 1,\\
Y_{j}&=V_1+\cdots +V_{j},\, j\geq 1,\\
Y_0&=0.
\end{align*}
First, denote $a(p) := \Ee_p(V_i)=pM=p[p^{-1}]$ and note that $a(p) \in [1/2, 1]$ for all $p \in (0,1)$.
Moreover, $Var(V_i) = \Ee_p(V_i^2)-(\Ee_p(V_i))^2= p(1-p)[p^{-1}]\leq 1$.
Now, consider the auxiliary process $\chi_n := Y_n - n(1-\kappa)/2 + c'$, with $\chi_0 = c'$.
Note that
$ \Ee_p(\chi_{n+1} - \chi_n \> | \> \chi_n=\chi) = a(p) - (1-\kappa)/2 > 0$. Moreover, if we define the stopping time $t_x = \min_{n\geq 0}\{\chi_n \leq x\}$, it follows from Theorem  2.5.18 in \cite{MPW17} that there exist $x_1$ and $\alpha > 0$
such that
\[ \Prob\left( \bigcap_{n=1}^{\infty}\{ Y_n \geq \frac{n}{2}(1-\kappa)-(c'-x_1)\} \right)= \Prob(t_{x_1} = \infty) \geq 1-\left(  \frac{1+x_1}{1 +\chi_0} \right)^{\alpha}. \]
So, for every $\eps \in (0,1)$ and $\kappa > 0$, we can appropriately choose $\alpha$ and $\chi_0 = c' > x_1$ such that the probability in the above display is greater than $\eps/2$. Analogously, if we define $\chi_n = -Y_n + n(1+\kappa)$, the upper bound can be found exactly as above, yielding

\begin{equation}
\label{pr1}
\mu_p\left(\bigcap_{n=1}^{\infty}\{\frac{n}{2}(1-\kappa)-c\leq Y_n\leq n(1+\kappa)+c\}\right)\geq \eps.
\end{equation}

Further, fix $n\geq M$. Let $m_n$ and $l_n$ be integers such that  $n=m_nM+l_n$, where $l_n<M$.
Then on event  $\bigcap_{n=1}^{\infty}\{\frac{n}{2}(1-\kappa)-c\leq Y_n\leq n(1+\kappa)+c\}$ the following bounds hold
\begin{equation}
\label{c1}
U_n\geq Y_{m_n}\geq \frac{1}{2}\left(\frac{n}{M}-\frac{l_n}{M}\right)(1-\kappa)-c\geq \frac{1}{2}np(1-\kappa)-c_1,
\end{equation}
and 
\begin{equation}
\label{c2}
U_n\leq Y_{m_n+1}\leq  \left(\frac{n}{M}+\frac{M-l_n}{M}\right)(1+\kappa)+c\leq np(1+\kappa)+c_2.
\end{equation} 
Inequalities (\ref{pr1}), (\ref{c1}) and (\ref{c2}) yield bound  (\ref{PL4_1}).

Remind that $M=[p^{-1}]$, and so,
$$p\sum_{i=0}^{M-1}e^{-\lambda U_i}\leq pM\leq 1.$$
By combining this bound with  bound~(\ref{PL4_1}), it follows that given $\eps \in (0,1)$ and $\kappa >0$ we can find $c_1>0$ such that 
with $\mu_p$-probability at least $\eps$
\begin{equation}
\label{C}
p\sum_{i=0}^{\infty} e^{-\lambda U_i}\leq 1+p\sum_{i=M}^{\infty} e^{-\lambda (\frac{1}{2}pi(1-\kappa)-c_1)}\leq C
\end{equation}
for some deterministic constant $C=C(\eps, \lambda)$ and all $p\in(0,1)$.
Therefore
$$\inf\limits_{p\in(0, 1)}\Ee_p\left(e^{-p\sum_{i=0}^{\infty}e^{-\lambda U_i}}\right)\geq \eps e^{-C}=\eps_1>0,$$
as required. 
\qed

\bigskip 

We are now ready to proceed with the proof of the lemma. Remind that  $\lambda_k=\lambda_{k+1}=:\lambda$.
\paragraph{{\it Proof of Part  1) of Lemma~\ref{L4}.}}
Recall that in this case  
 $\lambda_{k-1}<\lambda_{k}=\lambda_{k+1}=\lambda$, $\lambda\geq \lambda_{k+2}$
and  $\Gamma_k=\max_i\Gamma_i$.
Then,
\begin{equation}
\label{Sigma1}
\gamma_{k,1}Z(\zeta_1)=
\gamma_{k,1}\sum_{i=0}^{\infty} e^{\lambda_{k-1}(i-U_i)-\lambda i}\leq \sum_{i=0}^{\infty}e^{-(\lambda-\lambda_{k-1})i}\leq C_1<\infty,
\end{equation}
where $C_1>0$ is a deterministic constant and we used that $\gamma_{k, 1}\leq 1$.

Further, if $\lambda>\lambda_{k+2}$, then 
\begin{equation}
\label{Sigma21}
\gamma_{k,2}Z(\zeta_2)=\gamma_{k,2}\sum_{i=0}^{\infty}e^{\lambda_{k+2}U_i-\lambda i} \leq
\sum_{i=0}^{\infty}e^{-(\lambda-\lambda_{k+2})i}\leq C_2<\infty,
\end{equation}
where $C_2>0$ is a deterministic constant and we used that $\gamma_{k, 2}\leq 1$.
Then, using bounds (\ref{Sigma1}) and (\ref{Sigma21}) in lower bound   (\ref{lower1}) gives 
that  
$
\Px\left(A^{k, k+1}_{[1, \infty)}\right)\geq \eps
$
for some $\eps>0$, as claimed.

If $\lambda=\lambda_{k+2}$, then bound (\ref{Sigma21}) cannot be used,  and  we proceed as follows. 
Note  that in this case 
\begin{equation}
\label{Sigma22}
\gamma_{k,2}Z(\zeta_2)=
\gamma_{k,2}\sum_{i=0}^{\infty}e^{\lambda(U_i-i)} \leq
q\sum_{i=0}^{\infty}e^{\lambda(U_i-i)},
\end{equation}
as
\begin{equation}
\label{gammak2}
\gamma_{k, 2}=\frac{\Gamma_{k+2}}{\Gamma_{k}+\Gamma_{k+1}}\leq \frac{\Gamma_{k}}{\Gamma_{k}+\Gamma_{k+1}}=q=1-p,
\end{equation}
where $p$ is defined in~(\ref{p}).
Further, combining  bounds (\ref{Sigma1}) and (\ref{Sigma22}) in  (\ref{lower1})  we get that 
\begin{equation}
\label{eps_q}
\Px\left(A^{k, k+1}_{[1, \infty)}\right)\geq 
\eps_1 \Ee_p\left(e^{-q\sum_{i=0}^{\infty}e^{\lambda(U_i-i)}}\right)=\eps_1\Ee_p\left(e^{-p\sum_{i=0}^{\infty}e^{-\lambda U_i}}\right),
\end{equation}
where the  equality holds by symmetry. 
It is left to note that  the expectation in the right side of the last equation 
is  bounded below uniformly over $p\in(0, 1)$ by Part 2) of Proposition~\ref{PL4}.
\qed

\paragraph{{\it  Proof of Part 2) of Lemma~\ref{L4}.}}
Recall that in this case
$\lambda_{k-1}=\lambda_{k}=\lambda_{k+1}=\lambda\geq \lambda_{k+2}$,
$\Gamma_k=\max_i\Gamma_i$  and  $\Gamma_{k-1}\leq \Gamma_{k+1}$. 
These conditions give that 
$e^{\lambda_{k-1}(i-U_i)-\lambda i}=e^{-\lambda U_i}$, 
$e^{\lambda_{k+2}U_i-\lambda i}\leq e^{\lambda(U-i)}$, and 
\begin{equation}
\label{gammak1}
\gamma_{k, 1}=\frac{\Gamma_{k-1}}{\Gamma_k+\Gamma_{k+1}}\leq \frac{\Gamma_{k+1}}{\Gamma_k+\Gamma_{k+1}}=p.
\end{equation}
Recall also  that $\gamma_{k, 2}\leq q=1-p$ (see (\ref{gammak2})).
Using all these inequalities in lower  bound  (\ref{lower1})
gives  the following lower bound
\begin{equation}
\label{lower11}
\Px\left(A^{k, k+1}_{[1, \infty)}\right)
\geq \delta \Ee_p\left( e^{-p\sum_{i=0}^{\infty} e^{-\lambda U_i}}
e^{-q\sum_{i=0}^{\infty} e^{\lambda(U_i-i)}}  \right).
\end{equation}
We have already shown in  (\ref{C}) that for any $\eps\in(0, 1)$  there exists constant $C=C(\eps)>0$ such that 
\begin{equation}
\label{eUi<C}
\mu_p\left(p\sum_{i=0}^{\infty} e^{-\lambda U_i}\leq  C\right)\geq \eps \quad \text{and} \quad 
\mu_p\left(q\sum_{i=0}^{\infty} e^{\lambda (U_i-i)}\leq  C\right)\geq \eps
\end{equation}
for all $p$, where the second bound holds  by symmetry. 
Choosing $\eps>0.5$ we get that 
$$\mu_p\left(p\sum_{i=0}^{\infty} e^{-\lambda U_i}\leq  C,\, 
q\sum_{i=0}^{\infty} e^{\lambda (U_i-i)}\leq  C\right)\geq 2\eps-1>0,$$
for all $p$. 
Combining this bound with  equation (\ref{lower11}) we finally obtain that 
$
\Px\left(A^{k, k+1}_{[1, \infty)}\right) \geq \eps_2
$
for some $\eps_2>0$, as claimed. \qed

\subsubsection{Proof of Lemma~\ref{L5}}
\label{prL5}

\paragraph{{\it Proof of Part 1) of Lemma \ref{L5}}.}
Note that at every time a particle is added to site $k$ 
or $k+1$, the allocation rates at these sites are multiplied by $e^{\lambda}$.
In particular, if  a particle is added to site $k$, then the allocation rate at $k-1$ is multiplied by $e^{\lambda_{k-1}}$.
Otherwise, if a particle is added to site  $k+1$, then  the allocation rate at $k+2$ is multiplied by $e^{\lambda_{k+2}}$.
Other rates remain unchanged.
Thus, by  allocating a particle at $k$ or $k+1$, the sum of rates at $k, k+1$ and $k+2$ over the sum of rates at all other sites is increased by a multiple constant.
This yields the following exponential bound 
\begin{equation}
\label{expon}
\P_{\bx}\left(\bigcup_{i\neq k, k+1, k+2}A^i_{n+1}\bigg|A^{k, k+1}_{[1, n]}\right) \leq C_1e^{-C_2n},
\end{equation}
for some $C_1, C_2>0$.
In turn, bound (\ref{expon}) implies that with a positive probability (not depending on $\bx$) 
event $A_{[1,\infty)}^{k, k+1}\cup \{\tau_{k+2}<\infty\}$ occurs as claimed. 
Note also  that  events $A_{[1,\infty)}^{k, k+1}$ and $\{\tau_{k+2}<\infty\}$
are mutually exclusive.
Thus, with a positive probability either all particles will be allocated at $k$ and $k+1$, or a particle is eventually  placed at $k+2$.
Placing a particle at $k+2$ can violate  condition (\ref{A}) because 
the maximal allocating probability can be now attained  at sites  $k+2$ and $k+3$ as well. 
Part 1) of Lemma \ref{L5}  is proved. \qed

\paragraph{ {\it Proof of Part 2)  Lemma \ref{L5}.}}
Note that 
 $e^{\lambda_{k-1}(i-U_i)-\lambda i}<e^{-(\lambda-\lambda_{k-1})i}$ and 
 $\lambda_{k-1}<\lambda$. Consequently, for any $n$
\begin{equation}
\label{a1}
Z_n(\zeta_2)=\sum_{i=0}^n   e^{\lambda_{k-1}(i-U_i)-\lambda i}\leq \sum_{i=0}^{\infty} e^{-(\lambda-\lambda_{k-1})i}<C<\infty.
\end{equation}
Note also  that $\gamma_{k, 1}\leq 1$ and $\gamma_{k, 2}\leq 1$. Combining these inequalities with equation (\ref{a1})  and letting 
$n\to \infty$ in (\ref{lower_Z}) gives that
$$\Px\left(A^{k, k+1}_{[1, \infty]}\right)\geq \varepsilon_1\Ee_p\left(e^{-Z(\zeta_2)}\right),$$
for some $\eps_1>0$.
Further,   assumption $r<z_2$ implies that 
$\lambda_{k+2}p - \lambda < 0$.  Recall  that parameter $r=x_{k+2}-x_{k-1}$ takes integer values,
and $p=p(r)$ is a  monotonically increasing function of $r$. 
Let $r_{0}$ be the maximal integer such that $r<z_2$ and $p_0=p(r_0)$, so that $\lambda_{k+2}p_0-\lambda<0$.
  It follows from  Proposition \ref{P2} and 
Proposition \ref{zeta(p1)<zeta(p2)} that for all $0<p<p_0$
\begin{equation}
\label{p<p_0^-}
\Ee_{p}\left(e^{-Z(\zeta_2)}\right)\geq \Ee_{p_0}\left(e^{-Z(\zeta_2)}\right)>0,
\end{equation}
and, hence, 
$
\P_{\bx}\left(A_{[1, \infty]}^{k, k+1}\right)\geq \eps
$
for some uniform $\eps>0$ over configurations $\bx$ satisfying $r<z_2$.
Part 2) of Lemma \ref{L5} is  proved.\qed

\paragraph{
{\it Proof of Part 3) of  Lemma \ref{L5}}.}
We are going to use the following relaxation of upper bound (\ref{upper})
\begin{equation}
\label{34upper}
\Px\left(A^{k, k+1}_{[1, n+1]}\right)
\leq \Ee_p\left(
\prod\limits_{i=0}^n \frac{1}{1+
\gamma_{k,2}e^{\lambda_{k+2}U_i-\lambda i}}\right).
\end{equation}
Next, assumption $r\geq z_2$ implies that 
$\lambda_{k+2}p- \lambda \geq 0$. Therefore, by the strong law of large numbers,  we get  that $ \mu_p\text{-} a.s.$
$\lambda_{k+2}U_i-\lambda i\geq 0$  for infinitely many $i$
and, hence, 
$\prod_{i=0}^n \frac{1}{1+\gamma_{k,2}e^{\lambda_{k+2}U_i-\lambda i}}\to 0$. The product is bounded by $1$, therefore, 
by the Lebesgue's dominated convergence theorem,
 the expectation in the right side of (\ref{34upper}) tends to $0$ as $n\to \infty$, which implies 
 that 
\begin{equation}
\label{kk+1=0}
\Px\left(A_{[1, \infty]}^{k, k+1}\right)=\lim\limits_{n\to \infty}\Px\left(A_{[1, n+1]}^{k, k+1}\right)=0,
\end{equation}
as claimed. 
Note that  equation (\ref{kk+1=0}) combined with Part 1) of the lemma further yields  
that 
$\Px\left(\tau_{k+2}<w_{k}^+\right)>\epsilon$ for some $\epsilon$. \qed

\paragraph{{\it Proof of Part 4) of  Lemma \ref{L5}}.}
 Define 
\begin{equation}
\label{hatn}
\widehat n=\min\left(n:\gamma_{k,2}e^{\lambda_{k+2}S_n-\lambda n}\geq 1\right).
\end{equation}
In other words, $\hat n$ is the first time when the allocation rate at site $k+2$ exceeds the sum of allocation rates at sites $k$ and $k+1$, becoming therefore, the maximal rate.

Applying lower bound (\ref{Bern31_Z}) gives  that
\begin{equation}
\label{Bern32}
\Px\left(A^{k, k+1}_{[1, \widehat n]}\right)\geq \delta\Ee_p\left(
e^{- \gamma_{k,1}Z_{\widehat m -1}(\zeta_1)}
e^{- \gamma_{k,2}Z_{\widehat m -1}(\zeta_2)}\right),
\end{equation}
where 
\begin{equation}
\label{hatm}
\widehat m=\min\left(m:\gamma_{k, 2}e^{\lambda_{k+2}U_m-\lambda m}\geq 1\right).
\end{equation}
Equation (\ref{a1}) yields that  $ \gamma_{k,1}Z_{\widehat m -1}(\zeta_1)<Z(\zeta_1)<C<\infty$.
This allows us to rewrite bound (\ref{Bern32}) as follows 
$
\Px\left(A^{k, k+1}_{[1, \widehat n]}\right)\geq \varepsilon_2\Ee_p\left(e^{- \gamma_{k,2}Z_{\widehat m-1}(\zeta_2)}
\right),$
for some $\eps_2$.
By  assumption $r > z_2$. Let now $r_0$ be the minimal integer such that $r_0>z_2$ and $p_0=p(r_0)$.
Then $\lambda_{k+2}p- \lambda >\lambda_{k+2}p_0 - \lambda > 0$ for any  $p>p_0$.
It follows from Proposition~\ref{C3} and Proposition~\ref{zeta(p1)<zeta(p2)}
that for all $p>p_0$
\begin{equation}
\label{p>p_0^+}
\Ee_{p}\left(e^{-\gamma_{k,2} Z_{\whm-1}(\zeta_2)}\right)\geq \Ee_{p_0}\left(e^{-Z(\eta_2)}\right)>0,
\end{equation}
where $\eta_2$ is the sequence reciprocal to $\zeta_2$.
Hence 
$\Px\left(A^{k, k+1}_{[1,  \widehat n]}\right) \geq  \eps_2\epsilon_2>0$.

Next, recall event $B_k$ defined in (\ref{Bk}). 
Note that $A^{k, k+1}_{[1,  \widehat n]} \cap A_{\widehat n + 1}^{k+2}\subseteq B_k$, so that 
$\Px\left(B_k\right)\geq \eps_2\epsilon_2/N>0$ as well.

\bigskip 

It is left to show that the maximal rate $\max_i\Gamma_i$ relocates as described in  (\ref{l51bound}). 
Clearly, this is always the case if $\lambda<\min(\lambda_{k+2}, \lambda_{k+3})$. 
This might not be the case in the following particular situation. 
Namely, suppose that  $\lambda_{k+3} \leq \lambda$ and  initial configuration $\bx$ is such that $\Gamma_k = \max_i \Gamma_i$ and $\Gamma_{k+3}e^{\lambda_{k+3}} \geq \Gamma_k$. In this case, if $\tau_{k+2} = 1$, then the maximal rate might move   to $k+3$. 
However, note that $\tau_{k+2}\geq 2$ on event $A^{k, k+1}_{[1,  \widehat n]}$.
Indeed, by definition (\ref{hatn})  $\widehat n\geq 1$, and, hence, on this event $\tau_{k+2}\geq 2$ as $\tau_{k+2}>\widehat n$, so that at least one particle is deposited at $\{k, k+1\}$ by time $\tau_{k+2}$. It is not hard to check that placing one particle at $\{k, k+1\}$ makes impossible that relocation of $\max_i\Gamma_i$ to $k+3$ when $\lambda_{k+3}\leq \lambda$.\qed

\subsubsection{Proof of Lemma~\ref{L6}}
\label{prL6}
First, note that the proof of Part 1) of Lemma \ref{L6} is analogous to the proof of Part 1) of Lemma \ref{L5}  and we omit  technical details.  
For  simplicity of notation we denote $\lambda=\lambda_{k}=\lambda_{k+1}$ in the rest of the proof.
\paragraph{{\it Proof of Part 2) of Lemma \ref{L6}}.}
Recall  lower bound (\ref{lower1}) 
\begin{equation*}
\Px\left(A^{k, k+1}_{[1, \infty]}\right)\geq 
\delta \Ee_p\left( e^{-\gamma_{k,1}Z(\zeta_1)}
e^{- \gamma_{k,2}Z(\zeta_2)}\right).
\end{equation*}
Note that  $z_1 < r < z_2$ if and only if both  $\lambda_{k-1}(1-p)-\lambda<0$ and $\lambda_{k+2}p-\lambda<0$.
Therefore, it follows from Proposition \ref{P2} that $ \mu_p\text{-} a.s.$ both $Z(\zeta_1)<\infty$ and $Z(\zeta_2)<\infty$.
Consequently,
\[\Ee_p\left( e^{-\gamma_{k,1}Z(\zeta_1)}
e^{- \gamma_{k,2}Z(\zeta_2)}\right) \geq \Ee_p\left( e^{-Z(\zeta_1)}
e^{- Z(\zeta_2)}\right)\geq  \varepsilon(p)>0, \]
as $\gamma_{k, i}\leq 1$, $i=1,2$, so that $\Px\left(A^{k, k+1}_{[1, \infty]}\right)\geq \delta\varepsilon(p)$.
It is left to note that there is a finite number (depending only on $\lambda$'s) of possible values of integer-valued parameter  $r$ satisfying  $z_1<r<z_2$, and, hence, the same number of possible values of probability $p$.
Therefore, constant $\varepsilon(p)$ can be chosen as the minimal one for those values of $p$. 
This concludes the proof of the second part of the lemma.  \qed

\paragraph{{\it Proof of Part 3) of Lemma \ref{L6}}. } Let us start by noting the following. 
Assumption  $r \leq z_1$ implies that $\lambda_{k-1}(1-p)-\lambda\geq 0$, and
assumption  $r \geq z_2$ implies that  $\lambda_{k+2}p - \lambda  \geq 0$.
Therefore, the law of large numbers  yields that $\mu_p\text{-}a.s.$ at least one of the following events $\{\lambda_{k-1}(i-U_i)-\lambda i\geq 0\}$ and 
$\{\lambda_{k+2}U_i-\lambda i \geq 0\}$ occurs for infinitely many $i$.
Consequently, $ \mu_p\text{-} a.s.$  $\prod\limits_{i=0}^n \frac{1}{1+
\gamma_{k,1}e^{\lambda_{k-1}(i-U_i)-\lambda i}+
\gamma_{k,2}e^{\lambda_{k+2}U_i-\lambda i}}\to 0$,  as $n\to \infty$.
Using bound (\ref{upper}) and  the Lebesgue dominated convergence theorem,  we obtain that 
\[ \Px\left(A^{k, k+1}_{[1, n+1]}\right)
\leq \Ee_p\left(
\prod\limits_{i=0}^n \frac{1}{1+
\gamma_{k,1}e^{\lambda_{k-1}(i-U_i)-\lambda i}+
\gamma_{k,2}e^{\lambda_{k+2}U_i-\lambda i}}\right) \rightarrow 0,\]
as $n\to \infty$.
Hence,  
$\Px\left(A^{k, k+1}_{[1, \infty]}\right)=0$, and, hence, $\Px\left(D_k\right)\geq \varepsilon$, as claimed. \qed

\subsubsection{Proof of Lemma \ref{L7}}
\label{proofL7L8}

The proof here is similar to the proof of Part 4) of Lemma \ref{L5}.
The common starting point is the lower bound (\ref{Bern31_Z}) where $\tau$ and $\tilde \tau$ are appropriately chosen stopping times.

\paragraph{{\it Proof of Part 1) and 2) of Lemma \ref{L7}}. } First, note that the random variables $Z(\zeta_1)$ and $Z(\zeta_2)$ are finite if $\lambda_{k-1}(1-p)-\lambda<0$ and $\lambda_{k+2}p-\lambda<0$, respectively.
In fact, by our assumptions, precisely one of these conditions is necessarily satisfied so that one of $Z(\zeta_1)$ and $Z(\zeta_2)$ is almost surely finite.
Then we apply bound (\ref{Bern31_Z}) with the corresponding pair of stopping times  $(\tau, \tilde \tau)=(\whn_2, \whm_2)$ or $(\tau, \tilde \tau)=(\whn_1, \whm_1)$ respectively, where
\begin{align*}
\whn_1&=\min\left(n:\gamma_{k, 1}e^{\lambda_{k-1}(n-S_n)-\lambda n}\geq 1\right),\\
\whn_2&=\min\left(n:\gamma_{k, 2}e^{\lambda_{k+2}S_n-\lambda n}\geq 1\right),\\
\whm_{1}&=\min\left(m: \gamma_{k, 1}e^{\lambda_{k-1}(m-U_m)-\lambda m}\geq 1\right),\\
\whm_2&=\min\left(m: \gamma_{k, 2}e^{\lambda_{k+2}U_m-\lambda m}\geq 1\right).
\end{align*} 

For concreteness, consider the case where $\{k, k+1\}$ is of type 2 and $r>z_1\geq z_2$, in which case $\lambda_{k-1}(1-p)-\lambda<0$ and $\lambda_{k+2}p-\lambda>0$. 
Applying bound (\ref{Bern31_Z}) with $(\tau, \tilde \tau)=(\whn_2, \whm_2)$ yields that 
$$
\Px\left(A^{k, k+1}_{[1, \whn_2]}\right)\geq \delta \Ee_p\left( e^{- \gamma_{k, 1}Z_{\whm_2-1}(\zeta_1)}
e^{- \gamma_{k, 2}Z_{\whm_2-1}(\zeta_2)}\right).
$$
Condition $\lambda_{k-1}(1-p)-\lambda<0$  and  Proposition \ref{P2} imply that $Z(\zeta_1)<\infty \>\> \mu_p \text{-a.s.}$  
Therefore, we can bound $\gamma_{k, 1}Z_{\whm_2-1}(\zeta_1)\leq Z(\zeta_1)$, as $\gamma_{k, 1}\leq 1$.
Also, condition $\lambda_{k+2}p-\lambda>0$ and Proposition  \ref{P3} imply that $\gamma_{k, 2}Z_{\whm_2-1}(\zeta_2)<\infty \> \mu_p\text{-a.s.}$
Combining the above, we get to the following lower bound
$$
\Px\left(A^{k, k+1}_{[1, \whn_2]}\right)
\geq \delta \Ee_p\left( e^{-Z(\zeta_1)}
e^{- \gamma_{k, 2}Z_{\whm_2-1}(\zeta_2)}\right).
$$
Moreover, let $\eta_2$ be the sequence reciprocal to $\zeta_2$.
Then, applying  Proposition \ref{P3} again, we get that 
$Z(\eta_2)<\infty$ $\mu_p$-a.s., 
$ Z(\eta_2)\geq_{st}\gamma Z_{\whm-1}(\zeta_2)$ and
$$
\Ee_{p}\left( e^{-Z(\zeta_1)}
e^{- \gamma Z_{\whm-1}(\zeta_2)}\right)\geq \Ee_{p}\left( e^{-Z(\zeta_1)}
e^{- Z(\eta_2)}\right)>0.$$
Let us show that, when $r>z_1$, the expectation in the right side of the preceding display is  uniformly bounded below  over $p=p(r)$.
To this end, take the minimal integer $r_0$ such that $r_0>z_1$ so that condition $r>z_1$ implies $p>p_0=p(r_0)$, and, hence,
$\lambda_{k-1}(1-p)-\lambda<\lambda_{k-1}(1-p_0)-\lambda<0$ and 
$\lambda_{k+2}p-\lambda>\lambda_{k+2}p_0-\lambda>0$.
This implies the following.
First, consider the  random variable $Z(\zeta_1)$ with distribution determined by parameter $p_0$.
By Propositions~\ref{P2} and \ref{zeta(p1)<zeta(p2)}, it follows that $Z(\zeta_1)$ is almost surely finite, and, moreover, it stochastically dominates any other random variable $Z(\zeta_1)$ with distribution determined by $p>p_0$.
Second, consider the random variable $Z(\eta_2)$, where $\eta_2$ is a sequence reciprocal to sequence $\zeta_2$ whose distribution is determined by parameter $p_0$.
By Propositions~\ref{P2}, \ref{dominance} and \ref{zeta(p1)<zeta(p2)}, it follows that $Z(\eta_2)$ is almost surely finite and, moreover, it stochastically dominates any other random variable $Z(\eta_2)$, where $\eta_2$ is reciprocal to $\zeta_2$ whose distribution is determined by $p>p_0$.

Therefore,
$\Ee_{p}\left( e^{-Z(\zeta_1)}
e^{- Z(\eta_2)}\right)\geq  \Ee_{p_0}\left( e^{-Z(\zeta_1)}
e^{- Z(\eta_2)}\right)$.
Summarizing the above, we finally obtain that 
\begin{equation}
\label{Lemma82}
\Px\left(A^{k, k+1}_{[1, \whn_2]}\right)\geq 
\delta\Ee_{p_0}\left( e^{-Z(\zeta_1)}e^{-Z(\eta_2)}\right)>0.
\end{equation}
We have considered here only the case where $\{k,k+1\}$ is of type 2 and $r>z_1$, but by rearranging the stopping times above, one should note that for all the remaining cases stated in Parts 1) and 2) of Lemma \ref{L7}, the reasoning is exactly the same as above. \qed

\paragraph{{\it Proof of Part 3) of Lemma \ref{L7}}. }
Let us obtain  the lower bound in  Part 3) of Lemma \ref{L8}.
In this case  $\{k, k+1\}$ is a local minimum of type $2$ and 
$z_2<r<z_1$. The double inequality implies that 
both $\lambda_{k-1}(1-p)-\lambda>0$ and  $\lambda_{k+2}p-\lambda>0$. As a result, 
both $Z(\zeta_1)$ and  $Z(\zeta_2)$  are infinite.   
  In this case we modify  bound (\ref{Bern31_Z}) with stopping times $\tau=\widehat n=\min(\whn_1, \whn_2)$ and 
$\tilde \tau=\widehat m= \min(\whm_1, \whm_2),$  as follows
\begin{equation*}
\begin{split}
\Px\left(A^{k, k+1}_{[1, \widehat n]}\right)&\geq 
\delta \Ee_p\left( e^{- \gamma_{k,1}Z_{\widehat m-1}(\zeta_1)}
e^{- \gamma_{k,2}Z_{\widehat m-1}(\zeta_2)}  \right)\\
&\geq \delta \Ee_p\left( e^{-\gamma_{k,1}Z_{\widehat m_1-1}(\zeta_1)}
e^{-\gamma_{k,2}Z_{\widehat m_2-1}(\zeta_2)}  \right),
\end{split}
\end{equation*}
where in the last inequality we  bounded  $\whm=\min(\whm_1, \whm_2)$   by $\whm_1$ and $\whm_2$ respectively.
By Proposition \ref{P3} $\mu_p\text{-} a.s.$ both $ \gamma_{k,1}Z_{\widehat m_1-1}(\zeta_1)<\infty$ and 
$  \gamma_{k,2}Z_{\widehat{m}_2-1}(\zeta_2)<\infty$.
 Therefore,  $\Px\left(A^{k, k+1}_{[1, \widehat n]}\right)\geq \eps(p)>0$.   Further,  there are  finitely many integers 
$r$ such that $z_2<r<z_1$. Consequently, there are finitely many corresponding values of  probability  $p$, and
$\Px\left(A^{k, k+1}_{[1, \widehat n]}\right)\geq \varepsilon$ for some $\varepsilon>0$ uniformly over all values of $p$ in this finite set.  

Finally, relocation of the maximal rate in all cases covered by Lemmas~\ref{L7} can be shown by modifying the argument used in the proof of Part 4) of Lemma~\ref{L5}. \qed

\subsubsection{Proof of Lemma~\ref{L8}}
 We skip proofs of Parts 1) and 3) as they are analogous to the proofs of  Parts 1) and 3) 
of Lemma~\ref{L5}. 
Proofs of Parts 2) and 4) can be obtained by appropriately modifying proofs of  Parts 2) and 4) 
of Lemma~\ref{L5} and combining them with the  ideas in the proof  of Lemma~\ref{L4}. 
Modifications are due to condition $\lambda_{k-1}=\lambda$ implying that  $z_1=-\infty<z_2$ (see Remark~\ref{R1}).  

\paragraph{{\it Proof of Part 2) of Lemma~\ref{L8}.}}
Recall that in this case $r<z_2$, so that $\lambda_{k+2}p-\lambda<0$ and $p<p_0$, 
where $p_0$ is defined in Part 2) of Lemma~\ref{L5}.
 Repeating the proof of Part 2) of Lemma~\ref{L5} and using that $\gamma_{k, 1}\leq p$ and $\gamma_{k,2}\leq 1$ (see~(\ref{gammak1}) and~(\ref{gammak2})) we obtain  the following lower bound
\begin{equation}
\label{lower82}
\Px\left(A^{k, k+1}_{[1, \infty]}\right)\geq \Ee_p\left(e^{-pZ(\zeta_1)}e^{-Z(\zeta_2)}\right),
\end{equation}
Our assumptions imply that  both $Z(\zeta_1)$ and $Z(\zeta_2)$ are almost surely finite by Proposition \ref{P2}.
 Fix $\eps>0.5$, let  $C_1=C_1(\eps)>0$ be such that 
\begin{equation}
\label{eps_C}
\mu_p\left(pZ(\zeta_1)\leq  C_1\right)=\mu_p\left(p\sum_{i=0}^{\infty}e^{-\lambda U_i}\leq C_1\right)\geq \eps
\end{equation}
for all $p\in(0, 1)$ (see~(\ref{eUi<C})), 
and let  $C_2=C_2(\eps)$ be  such that $\mu_{p_0}\left(Z(\zeta_2)\leq  C_2\right)\geq \eps$.
The last inequality yields that  
$\mu_{p}\left(Z(\zeta_2)\leq  C_2\right)\geq \mu_{p_0}\left(Z(\zeta_2)\leq  C_2\right)\geq \eps,$
as $Z(\zeta_2)$, with distribution determined by parameter $p_0$, dominates any random variable 
$Z(\zeta_2)$ with distribution determined by parameter $p<p_0$.
Finally,  by using the same elementary argument as in the proof of Lemma~\ref{L4},  we get  that 
$\mu_p\left(pZ(\zeta_1)\leq  C_1, Z(\zeta_2)\leq C_2\right)\geq 2\eps-1$,
 which implies  that the expectation in the right side 
of (\ref{lower82}) is bounded below away from zero, so
 that 
$
\P_{\bx}\left(A_{[1, \infty]}^{k, k+1}\right)\geq \eps_1
$
for some uniform $\eps_1>0$ over configurations $\bx$ satisfying $r<z_2$. \qed

\paragraph{{\it Proof of Part 4) of Lemma~\ref{L8}.}}
Recall that in this case  $r>z_2$, so that $\lambda_{k+2}p-\lambda>0$ and $p>p_0$, 
where $p_0$ is now  defined in Part 4) of Lemma~\ref{L5}.
 Repeating the proof of Part 4) of Lemma~\ref{L5} and using again  that $\gamma_{k, 1}\leq p$
we obtain  the following lower bound
$$
\Px\left(A^{k, k+1}_{[1, \widehat n]}\right)\geq \delta\Ee_p\left(
e^{-pZ(\zeta_1)}e^{- \gamma_{k,2}Z_{\widehat m -1}(\zeta_2)}\right),
$$
where 
$ \widehat n$ and  $\widehat m$ are defined in (\ref{hatn}) and (\ref{hatm}) respectively.
Our assumptions imply that  both $Z(\zeta_1)$ and $Z_{\widehat m -1}(\zeta_2)$ are almost surely finite by Propositions \ref{P2} 
and~\ref{P3}. 
Further, Proposition~\ref{P3} yields that 
\begin{equation}
\label{lower84}
\Px\left(A^{k, k+1}_{[1, \widehat n]}\right)\geq \delta\Ee_p\left(
e^{-pZ(\zeta_1)}e^{-Z(\eta_2)}\right),
\end{equation}
where $\eta_2$ is the  random sequence reciprocal to $\zeta_2$.

Let  $\eps>0.5$ and   $C_1=C_1(\eps)>0$ be such that  (\ref{eps_C}) holds, 
and  let  $C_2=C_2(\eps)$ be  such that 
$\mu_{p_0}\left(Z(\eta_2)\leq  C_2\right)\geq \eps.$
The last inequality yields that  
$$\mu_{p}\left(Z(\eta_2)\leq  C_2\right)\geq \mu_{p_0}\left(Z(\eta_2)\leq  C_2\right)\geq \eps,$$
as $Z(\eta_2)$, with distribution determined by parameter $p_0$, dominates any random variable 
$Z(\eta_2)$ with distribution determined by parameter $p>p_0$.

As at the  same stage of the proof in  Part 2) we can now  conclude that the expectation 
in the right side of (\ref{lower84}) is bounded below away from zero, which implies  that 
$
\P_{\bx}\left(A_{[1, \widehat n]}^{k, k+1}\right)\geq \eps_2
$
for some uniform $\eps_2>0$ over configurations $\bx$ satisfying $r>z_2$. \qed

\subsubsection{Proof of Corollary~\ref{critical}}

The critical cases where $r = z_1$ or $r = z_2$ need to be treated separately since these cases can not be proven directly by the above arguments. However, by a slight modification one can amend the proof of each lemma in order to encompass such critical cases.

The modification is the same for all lemmas, but for the sake of concreteness let us consider the critical case described in Part 3) of Lemma~\ref{L5} assuming that $r=z_2$.  
We start by commenting on the same effect that we already discussed in the proof of Part 4) of Lemma~\ref{L5}. Namely, recall that if
$\lambda_{k+3}<\lambda_k=\lambda_{k+1}$, $\Gamma_{k+3}e^{\lambda_{k+3}}\geq \Gamma_k$, and  $\Gamma_k=\max_i\Gamma_i$, then 
$\tau_{k+2} = 1$ makes the maximal rate  move  to $k+3$. 
One can check that the above situation is the only one that can possibly relocate the maximal rate to a site with smaller $\lambda$. 
In order to avoid such case, it is simply a matter of placing a particle at $k$ at the first step, which can be done with probability at least $1/N$.
Therefore, without loss of generality we can exclude this case.

Next, if at time $\tau_{k+2}$ the maximal rate relocates either to $k+2$, or to $k+3$ (provided $\lambda_{k+3}>\lambda_k=\lambda_{k+1}$) then we are done.
Suppose the opposite, namely, that at time $\tau_{k+2}$ the maximal allocation rate remains where it was, that is, at $k$ or at $k+1$.
It is left to note  that given event $A_{[1, \tau_{k+2}-1]}^{k, k+1}$, placing a particle at site $k+2$ at moment $\tau_{k+2}$ increases the configuration parameter $r=x_{k+2}-x_{k-1}$ by $1$, so that the resulting  configuration is such that $r>z_2$.
By Part 4) of Lemma~\ref{L5}, the next allocated particles at $\{k, k+1\}$ will end up by relocating the maximal rate as prescribed.

Other critical cases can be handled similarly, and we skip straightforward technical details.

\section{Proof of Theorem~\ref{T1}}
\label{proofT1}

The idea of the proof goes briefly as follows. Given any initial state $X(0) = \bx$, the site $k$ where $\Gamma_k(\bx)=\max_{i=1,\ldots, N}(\Gamma_i(\bx))$ is identified. Then, a particle allocation strategy is drawn so that it always results in localization of growth as described in Theorem \ref{T1}.
Lemmas~\ref{L1}-\ref{L8} enable us to identify the corresponding strategy for each particular case and bound its probability from below \emph{uniformly over initial configurations}.
Should a particular strategy fail to happen, which means that at a certain step $n$ a particle is not allocated according to that strategy, but somewhere else, a new one is drawn and this procedure reiterates from $X(n)$.
Since there is a finite number of possible strategies it follows from the renewal argument below that almost surely one of them eventually succeeds. 

In what follows, when referring to  Lemma~\ref{L2} or one of Lemmas~\ref{L4}-\ref{L8}, this automatically includes the symmetric cases by re-labelling the graph in reverse order (as explained in Remark~\ref{R2}). 
Also, local minima of size $2$ and type $1$ automatically include the limiting case described in Remark~\ref{R1}.

\medskip 
Let $X(n) = \bx$ be a fixed and arbitrary configuration, and:

\medskip

1) Assume that   
$
\Gamma_k(\bx)=\max_{i=1,\ldots, N}(\Gamma_i(\bx))$ and $\lambda_{k-1}\neq \lambda_k \neq \lambda_{k+1}$.

\indent \indent
1.1) Let  $k$ be  a local maximum. By Lemma \ref{L1}, with positive probability, all subsequent particles are allocated at $k$. 

\indent \indent
1.2)
 Let $k$ be either  a growth point,  or a local minimum. 
By Lemmas~\ref{L2} and \ref{L3}, with positive probability, the maximal rate relocates in finite time to one of its nearest neighbours having parameter $\lambda > \lambda_k$.

\medskip 

2)
Assume that $\Gamma_k(\bx)=\max_{i=1,\ldots, N}(\Gamma_i(\bx))$ and that additional assumptions of Lemma~\ref{L4} are satisfied.
Lemma~\ref{L4} yields that, with positive probability, all subsequent particles are allocated at sites $\{k, k+1\}$

\medskip 

3)
Assume that $\max(\Gamma_{k}(\bx), \Gamma_{k+1}(\bx))=\max_{i}\Gamma_i(\bx)$, where  $\{k, k+1\}$ is either a saddle point, or a local minimum of size $2$  and type $1$.
Additional assumptions on $\bx$, as described in Part 2) of  Lemmas~\ref{L5}, \ref{L6} and \ref{L8}, guarantee
 that, with positive probability, all subsequent particles are allocated at sites $\{k, k+1\}$.

\medskip 

4)
Assume that $\max(\Gamma_{k}(\bx), \Gamma_{k+1}(\bx))=\max_{i}\Gamma_i(\bx)$, where  $\{k, k+1\}$ is either a saddle point of size $2$, or a local minimum of size $2$ of {\it either} type.
Assume also that configuration $\bx$ is such that assumptions as in the preceding item do \emph{not} hold.
Such cases are covered by Lemmas:
~\ref{L5}, Part 3) and 4);
~\ref{L6} Part 3); 
~\ref{L7} and \ref{L8};
and finally, \ref{L8} Part 3) and 4) complemented by Corollary~\ref{critical}.
In all those cases, with positive probability, the maximal rate eventually relocates in a random but {\it finite} time to a site with larger parameter $\lambda$. 

\medskip

5)
Finally, for the remaining cases of local minima, maxima or saddle points of size greater than $2$, it is not hard to check that such cases can be reduced to one, or a combination, of the above items. 

\medskip 

Thus, for every configuration $\bx$ and every set of positive real parameters $\Lambda = (\lambda_k)_{k=1}^N$, we have identified two types of events.
First, there are events resulting in localisation of growth at either a single site or a pair of neighbouring sites (as described in Theorem~\ref{T1} Part 1) and 2) respectively).
Call such events $\L$-events.
Second, there are events resulting in relocation of the maximal rate. 
Call such events $\Rr$-events.

The next step of the proof is to define a sequence of random moments of time $(T_j)_{j \geq 0}$ called renewal moments. 
First, set  $T_0=0$.
Now, given $T_j$, let us define $T_{j+1}$. 
Suppose that at time $T_j$ the process is at state $\bx$.
We identify an event $R_1\ldots R_mL$ (strategy) formed by a sequence of $m$ $\Rr$-events (possibly none) ending at an $\L$-event.
At the fist moment of time $t>T_j$ a particle is not allocated according to $R_1\ldots R_mL$, we set $T_{j+1}=t$.

Note that $\Rr$-events are defined in a way so that the maximal rate always relocates to a site with strictly larger parameter $\lambda$.
It follows that the number of $\Rr$-events preceding any $\L$-event is bounded by the number of different values of $\lambda_i, \> i=1\ldots N$.
Then, by Lemmas~\ref{L1}-\ref{L8}, probabilities of events $R_1\ldots R_mL$ are bounded below uniformly over configurations, where $m \leq N$.

Further, let $j_{max}:= \max\{ j \geq 0 : T_j < \infty\}$.
Lemmas~\ref{L1}-\ref{L8} imply the existence of an uniform bound $\epsilon>0$ such that 
$\P(T_j = \infty)\geq \epsilon$ on $\{T_{j-1} < \infty\}$.
Therefore, $\P(T_j < \infty) \leq 1 - \epsilon$ on $\{T_{j-1} < \infty\}$, or equivalently, $\P(j_{max} \geq j \,|\, j_{max} \geq j-1) < 1-\epsilon$.
Thus,  $\P(j_{max}< \infty) = 1$.
This implies that $T_j = \infty$ for some $j$, so that,  with probability one,
 a certain allocation strategy $R_1\ldots R_mL$  eventually succeeds, that is 
the growth  process localises as claimed.

Finally, the long term behaviour of ratio $X_{k+1}(n)/X_k(n)$ described in item  ii) of the theorem
is implied by  the law of large numbers 
for the Binomial distribution.
This follows  straightforwardly from the proofs of Lemma~\ref{L4} and Parts 2) of Lemmas~\ref{L5}, \ref{L6} and \ref{L8}. 
The theorem is proved. \qed

\section*{Acknowledgements}
M.C.'s research is supported by CNPq [grant number 248679/2013-9]. \\
M.V.'s research is supported by CNPq [grants number 305369/2016-4 and\\ 302593/2013-6]. \\
V.S.'s research is partially supported by LMS [grant number 41562].\\
The authors are grateful to  Andrew Wade for suggestions that  improved the presentation.

\bibliography{GrowthPaperBibliography} 
\bibliographystyle{acm}
\end{document}